\documentclass[reqno,titlepage]{amsart}
\usepackage[english]{babel}
\usepackage{amsmath,amssymb,amsthm,amsaddr}
\usepackage{latexsym,mathtools,graphicx,epstopdf,hyperref}
\usepackage{enumerate}
\newlength\figureheight
\newlength\figurewidth
\mathtoolsset{showonlyrefs}
\usepackage{tikz,pgfplots}
\pgfplotsset{compat=newest}
\pgfplotsset{plot coordinates/math parser=false}

\begin{document}
\newcommand{\C}{{\mathbb{C}}}
\newcommand{\LOelr}{{\rL^2\left(\Omega_{r_\mathrm{e}}\right)}}
\newcommand{\LOe}{{\rL^2(\Omega_{r_\mathrm{e}})}}
\newcommand{\LOplr}{{\rL^2\left(\Omega_{r_\mathrm{p}}\right)}}
\newcommand{\LOp}{{\rL^2(\Omega_{r_\mathrm{p}})}}
\newcommand{\LOqlr}{{\rL^2\left(\Omega_{r_\mathrm{q}}\right)}}
\newcommand{\LOq}{{\rL^2(\Omega_{r_\mathrm{q}})}}
\newcommand{\LOslr}{{\rL^2\left(\Omega_{r_\mathrm{s}}\right)}}
\newcommand{\LOs}{{\rL^2(\Omega_{r_\mathrm{s}})}}
\newcommand{\LR}{{\rL^2(R)}}
\newcommand{\N}{{\mathbb{N}}}
\newcommand{\Op}{{\Omega_{r_\mathrm{p}}}}
\newcommand{\Oq}{{\Omega_{r_\mathrm{q}}}}
\newcommand{\Os}{{\Omega_{r_\mathrm{s}}}}
\newcommand{\R}{{\mathbb{R}}}
\newcommand{\Sz}{{\mathbb{S}^2}}
\newcommand{\cI}{{\mathcal{I}}}
\newcommand{\cL}{{\mathcal{L}}}
\newcommand{\cP}{{\mathcal{P}}}
\newcommand{\cR}{{\mathcal{R}}}
\newcommand{\cS}{{\mathcal{S}}}
\newcommand{\cT}{{\mathcal{T}}}
\newcommand{\cX}{{\mathcal{X}}}
\newcommand{\cY}{{\mathcal{Y}}}
\newcommand{\cZ}{{\mathcal{Z}}}
\newcommand{\doe}{{\,\mathrm{d}\omega(\eta)}}
\newcommand{\dox}{{\,\mathrm{d}\omega(\xi)}}
\newcommand{\doz}{{\,\mathrm{d}\omega(\zeta)}}
\newcommand{\dx}{{\,\mathrm{d}x}}
\newcommand{\dy}{{\,\mathrm{d}y}}
\newcommand{\dz}{{\,\mathrm{d}z}}
\newcommand{\rL}{{\mathrm{L}}}
\newcommand{\sumlmL}{{\sum_{l=0}^L\sum_{m=-l}^l}}
\newcommand{\sumlm}{{\sum_{l=0}^\infty\sum_{m=-l}^l}}
\newcommand{\supp}{{\mathop{\mathrm{supp}}\,}}
\newcommand{\xip}[2]{{\left\langle #1, #2 \right\rangle_\cX}}  
\newcommand{\xoip}[2]{{\left\langle #1, #2 \right\rangle_{\cX_1}}}  
\newcommand{\xtip}[2]{{\left\langle #1, #2 \right\rangle_{\cX_2}}}  
\newcommand{\yip}[2]{{\left\langle #1, #2 \right\rangle_\cY}}
\newcommand{\yoip}[2]{{\left\langle #1, #2 \right\rangle_{\cY_1}}}
\newcommand{\ytip}[2]{{\left\langle #1, #2 \right\rangle_{\cY_2}}}
\newcommand{\zip}[2]{{\left\langle #1, #2 \right\rangle_\cZ}}
\theoremstyle{plain}
\newtheorem{theorem}{Theorem}[section]
\newtheorem{lemma}[theorem]{Lemma}
\newtheorem{cor}[theorem]{Corollary}
\theoremstyle{definition}
\newtheorem{defin}[theorem]{Definition}
\newtheorem{example}[theorem]{Example}
\theoremstyle{remark}
\newtheorem{rem}[theorem]{Remark}

\title[Regularizing Inverse Problems with Slepian Wavelets]
{A General Approach to Regularizing Inverse Problems with Regional Data using Slepian Wavelets}

\author{Volker Michel}
\address{Geomathematics Group\\
Department of Mathematics\\
University of Siegen\\
Germany}
\email[V.\ Michel]{michel@mathematik.uni-siegen.de}
\urladdr{www.geomathematics-siegen.de}

\author{Frederik J.\ Simons}
\address{Department of Geosciences\\
Guyot Hall\\
Princeton University\\
U.S.A.}
\email[F.J.\ Simons]{fjsimons@alum.mit.edu}
\urladdr{www.frederik.net}

\begin{abstract}
Slepian functions are orthogonal function systems that live on
subdomains (for example, geographical regions on the Earth's surface,
or bandlimited portions of the entire spectrum). They have been firmly
established as a useful tool for the synthesis and analysis of
localized (concentrated or confined) signals, and for the modeling and
inversion of noise-contaminated data that are only regionally
available or only of regional interest. In this paper, we consider a
general abstract setup for inverse problems represented by a linear
and compact operator between Hilbert spaces with a known
singular-value decomposition (svd). In practice, such an svd is often
only given for the case of a global expansion of the data (e.g.\ on
the whole sphere) but not for regional data distributions. We show
that, in either case, Slepian functions (associated to an arbitrarily
prescribed region and the given compact operator) can be determined
and applied to construct a regularization for the ill-posed regional
inverse problem. Moreover, we describe an algorithm for constructing
the Slepian basis via an algebraic eigenvalue problem. The obtained
Slepian functions can be used to derive an svd for the combination of
the regionalizing projection and the compact operator. As a result,
standard regularization techniques relying on a known svd become
applicable also to those inverse problems where the data are
regionally given only. In particular, wavelet-based multiscale
techniques can be used. An example for the latter case is elaborated
theoretically and tested on two synthetic numerical examples.
\end{abstract}


\maketitle

\section{Introduction}\label{Sec:Intro}
In a wide range of scientific applications, concentrated in but not
confined to the geosciences, regional modelling from global data has
become increasingly important, for a variety of reasons. For example,
regional phenomena like the melting of the Greenland or Antarctica ice
sheets are being studied on the basis of global satellite
(potential-field, e.g.\ gravity) data \cite{Kusche2015}.
Alternatively, geophysical data could be of regionally varying
quality, either in terms of their measurement density, or owing to
spatial variations in signal-to-noise ratios. Finally, localization
and regionalization may be part of a strategy to \lq divide and
conquer\rq\ data domains, which is often a necessity for solving the
kinds of problems that involve the large data volumes with which the
geosciences are routinely confronted.

In this general context \cite{FreMicSim2017}, the use of localized
trial functions has proven to be useful. One among the many ways by
which such Ansatz functions can be constructed, the idea behind the
``Slepian'' approach is as follows.  Taking $R$ to be a subdomain, a
portion of a complete domain~$D$ (e.g.\ an interval on the set of real
numbers, or a spherical cap on the surface of a ball), we determine
the function~$F$ that maximizes the fraction
\begin{equation}\label{Eq:energyratio}
\lambda\coloneqq\frac{\int_R [F(x)]^2\dx}{\int_D [F(x)]^2\dx}
,
\end{equation}
the quotient of the squared $\rL^2$-norms of $F$ on~$R$ and
on~$D$. For practical purposes, the choice of $F$ is restricted to a
finite-dimensional space. This is achieved, for example, by assuming a
band\-limit for $F$. The first notions of Slepian functions treated
the case of the real line, and appeared in the literature in the early
1960s, in the work
by \cite{LandauPollak1961,Slepian1964,SlepianPollak1961}, who were
concerned with problems in communication theory. In the late 1990s,
Slepian functions on the sphere were derived for use in geodesy and
planetary
science \cite{AlbSanSne1999,DahlenSimons08,Simons09,SimonsDahlen06,SimDahWie06,WieczorekSimons05,WieczorekSimons07}. In
parallel, a few alternative approaches, using different measures of
optimality, have been developed for constructing approximating
structures on the sphere, see, for
example, \cite{Fernandez2007,FernandezPrestin2003,Michel2011}.

If $D$ is the 2-sphere, the functions $F$ can be expanded in the
well-known $\rL^2(D)$-orthonormal system of spherical harmonics
$Y_{l,m}$ of degree $l$ and order $m$ (see
e.g.\ \cite{OxfordBuch,Michelbook2013,Mueller1966}) up to a fixed
maximal degree~$L$,
\begin{equation}\label{Eq:expansion}
F(\xi)=\sumlmL \langle F, Y_{l,m}\rangle_{\rL^2} Y_{l,m}(\xi),
\quad |\xi|=1\,.
\end{equation}
The maximization problem~\eqref{Eq:energyratio} leads to an algebraic
eigenvalue problem, whose eigenvectors are vectors with the expansion
coefficients of~$F$ in the chosen basis (in the above case, the
$Y_{l,m}$), and whose eigenvalues are the ratios~$\lambda$
in~\eqref{Eq:energyratio}. Since the corresponding matrix is Gramian
and, therefore, symmetric, an orthonormal basis of eigenvectors
spanning the entire space of possible expansion coefficient vectors
can be found. Owing to Parseval's identity, the functions that
correspond to the expansion coefficient vectors also constitute an
orthonormal basis for the (bandlimited) space of considered
functions. As a consequence, the previously used basis can be replaced
by a new basis, the \lq Slepian\rq\ basis, whose elements are sorted
according to their localization~$\lambda$ over the subdomain~$R$. This
new basis is also orthogonal in the sense of $\rL^2(R)$, which
simplifies the expansion of bandlimited signals that are restricted to
the subdomain~$R$.

Recently, Slepian functions have revealed themselves to be also useful
for the regularization of inverse problems in geophysics. For
example, \cite{PlattnerSimons2017} addressed the downward continuation
of a gravity or magnetic field from regionally given gradients of the
potential at satellite altitude. Using $\cT F=G$ to represent the
inverse problem, involving an operator~$\cT $, a given function~$G$,
and an unknown function~$F$, \cite{PlattnerSimons2017} constructed
Slepian basis functions via the maximization of
\begin{equation}\label{Eq:enratinvprobl}
\tilde{\lambda}\coloneqq \frac{\int_R [\cT F(x)]^2 \dx}{\int_D [F(x)]^2\dx}\,.
\end{equation}
Here, $R$ need not be a subset of~$D$ any more, but, rather, is the
domain of functions in the range of~$\cT $. In the particular case
considered by \cite{PlattnerSimons2017}, $R$ is a region at satellite
altitude where data are being collected, and $D$ represents the
(spherical) Earth's surface.

In this paper, we will show that eqs.~\eqref{Eq:energyratio}
and \eqref{Eq:enratinvprobl} can be seen as particular examples of a
more general approach to the construction of Slepian functions for
inverse problems. In particular, in the typical application scenario,
one has an inverse problem $\cT F=G$ for which an svd is known if and
when~$G$ is given on a domain $\tilde{D}$. When $G$ is only given on a
subdomain $R\subset \tilde{D}$, we show that Slepian functions can be
used to derive an svd also for the restricted case $\cP\cT
F=\left. G\right|_R$, with a corresponding projection operator
$\cP$. The knowledge of such an svd opens the door to various
established regularization methods.

To the knowledge of the authors, there are only a few other
publications which use Slepian functions for inverse problems. For
example, another approach which addresses the singular-value
decomposition of the operator is developed in \cite{Khare2007} for
functions on the real line and a particular integral
operator. Moreover, in \cite{AbMoPh2015}, the gravitational potential
is expanded in spherical Slepian functions. The result is used as the
given right-hand side for an inverse problem, where point masses are
reconstructed which approximately generate the corresponding regional
gravitational potential. Examples in other application domains
are \cite{Etemadfard+2016,Jahn+2013,Mitra+2006,Sharifi+2014}.

Other systems of localized trial functions have been used for inverse
problems as well. This includes, in particular, wavelet
methods \cite{Gerhards2014a}. It would be beyond the scope and size of
this article to give a complete survey of such papers here. Examples
of other works where wavelets have been used for inverse problems on
the sphere
are \cite{HielscherSchaeben2008,Simonsetal2011,StaMurFad2015,Vareschi2014}. In
\cite{FreGloLit1999,MichelHab} it was shown that a wavelet-based regularization
can be constructed if the svd of the forward operator is known. 
Therefore, we use these
latter papers as a motivation for establishing a Slepian-based wavelet
method for inverse problems with regional data.

The outline of this paper is as follows: in Section~\ref{Sec:Not}, we
introduce some basic notation. The general setup of a linear compact
operator between two Hilbert spaces is described in
Section~\ref{Sec:Setting}. For this scenario, we explain the
construction of Slepian functions in Section~\ref{Sec:Slepian}. Since
the general setting includes also infinite-dimensional spaces but
numerical implementations are only possible for finite dimensions, the
practical specifics are discussed in Section~\ref{Sec:findim}. Since
the setting of \cite{PlattnerSimons2017} also includes an inverse
problem where data originating from two different kinds of sources are
being inverted, we show in Section~\ref{Sec:coup} how such coupled
problems can be integrated into the general scenario. In
Section~\ref{Sec:solv_ip} we describe an algorithm for determining the
Slepian functions and calculating the svd of the restricted
(projected) forward operator. Motivated by some known results for
Slepian functions on particular domains, we show in
Section~\ref{Sec:Fredholm} how Slepian functions can be used to
establish Fredholm integral operators for the forward and the inverse
operator. In particular, we also show how scaling functions and
wavelets can be constructed from the Slepian functions, and we prove
convergence and stability of the method. This multiscale
regularization technique is then applied to two inverse problems and
tested numerically for synthetic data sets in
Section~\ref{Sec:numerics}. Finally, in Section~\ref{Sec:concl}, we
offer conclusions and an outlook on future research.

\section{Notation}\label{Sec:Not}
As usual, $\N$ represents the set of all positive integers, where
$\N_0\coloneqq\N\cup\{0\}$, and $\R$ and $\C$ stand for the fields of
all real and complex numbers, respectively. A 2-sphere with radius
$r>0$ in $\R^3$ and centre~$0$ is denoted
\begin{equation}
\Omega_r \coloneqq \left\{\xi\in\R^3 \,\middle|\, |\xi| =r \right\}\,.
\end{equation}
We write $\Omega\coloneqq \Omega_1$ for the unit sphere, $r=1$.
Moreover, if $D\subset\R^n$ is measurable, then $\rL^2(D)$ is the
Hilbert space of square-integrable functions, where almost everywhere
equal functions are collected in equivalence classes.

\section{Setting}\label{Sec:Setting}
As we mentioned in the Introduction, we will present a general setup
for Slepian functions. For this purpose, we introduce here an abstract
setting which will serve as a starting point.  We have three
non-trivial Hilbert spaces $(\cX,\xip{\cdot}{\cdot})$,
$(\cY,\yip{\cdot}{\cdot})$, and $(\cZ,\zip{\cdot}{\cdot})$, with the
following additional assumptions.

\begin{itemize}
  \item There exists an isometric embedding (an injection)
  $\iota:\cZ\hookrightarrow\cY$, i.e.  
  \begin{equation} \yip{\iota\left(F_1\right)}{\iota\left(F_2\right)}
  = \zip{F_1}{F_2} \quad \text{for all }
  F_1,F_2\in\cZ\,.\label{Eq:isometry} 
  \end{equation} 
  We, therefore, consider $\cZ$ to be a subset of $\cY$ by associating
  $\cZ$ with $\iota(\cZ)$. Since $\iota$ is isometric and
  $(\cZ,\zip{\cdot}{\cdot})$ is a Hilbert space, also
  $(\iota(\cZ),\yip{\cdot}{\cdot})$ is a Hilbert space, namely, a
  Hilbert subspace of $(\cY,\yip{\cdot}{\cdot})$.
  \item There exists a projection $\cP:\cY\to\cZ$, in the sense that
  (with \lq $\cZ\subset\cY$\rq)
  \begin{equation}
   \cP(\cP G) = \cP G \quad \text{for all } G\in\cY, 
  \end{equation} 
  such that $\cP\circ\iota =\mathrm{Id}_\cZ$, in other words, $\cP$ inverts the
  embedding.
\end{itemize}
For a better understanding, we discuss an example of an application.

\begin{example}\label{Ex:upcontspaces}

A typical challenging inverse problem in the geosciences is the
downward continuation problem (see
e.g.\ \cite{Schachtschneider+2010,Schneider97,TelschowDiss}). As considered
by \cite{PlattnerSimons2017}, a harmonic potential (e.g.\ the
gravitational or magnetic potential) is given on a sphere with
radius~$r_\mathrm{s}$ (e.g.\ the satellite orbit), and the task is to
determine the potential on the surface of the planet (a sphere with
radius $r_\mathrm{p}$). In this case, we might choose
\begin{equation}
\cX=\LOplr,  \quad \cY=\LOslr, \quad \cZ=\LR,
\end{equation}
as spaces where $R\subset\Os$ is a subdomain, also a 2-dimensional
surface. For example, $R$ could be an area of limited access by
measurement, or to which the analysis of the potential is
restricted. The canonical embedding would then be
$\iota: \LR\hookrightarrow\LOs$ with
\begin{equation}
\left[\iota (F)\right] (x) \coloneqq \begin{cases}
F(x), & x\in R\\
0, & x\not\in R
\end{cases},
\qquad x\in \Os.
\end{equation}
It is clear that, for real $F_1,F_2\in\LR$, we have
\begin{align}
\left\langle \iota (F_1),\iota (F_2) \right\rangle_\LOslr 
&= \int_\Os \left[\iota (F_1)\right](\xi)\,  \left[\iota (F_2)\right](\xi)\dox\\
&= \int_R F_1(\xi) \, F_2(\xi) \dox \\
&= \left\langle F_1, F_2 \right\rangle_\LR\, .
\end{align}
The projection $\cP\colon\LOs\to\LR$ would simply be the restriction
\begin{equation}
\cP\colon G\mapsto \left. G\right|_R\, .
\end{equation}
It is similar to the restriction operator used in \cite[their Eq.\
(4.22)]{SimDahWie06}.
\end{example}

Let us return to the general setting again.
\begin{lemma}
We have
\begin{equation}
\left.\iota\circ\cP\right|_{\iota(\cZ)} = \mathrm{Id}_{\iota(\cZ)}
,\end{equation}
and $\iota\circ\cP$ is a projection onto $\iota(\cZ)$.
\end{lemma}
This lemma easily follows from the required properties above.


We will now continue with the abstract setting for the inverse
problem. For this purpose, we also assume that we have a compact
operator $\cT\colon\cX\to\cY$ with a known svd
\begin{equation}
\cT F = \sum_n \sigma_n \xip{F}{u_n} \,v_n, \quad F\in\cX,
\label{Eq:svd}
\end{equation}
where $(\sigma_n)_n\subset\C$ satisfies $\sigma_n\neq 0$ for all
$n$. Moreover, as usual for an svd, $(u_n)_n$ and $(v_n)_n$ are
orthonormal systems in $\cX$ and $\cY$, respectively.

Furthermore, we have an inverse problem $\cT F=\tilde{G}$, where
$\tilde{G}\in\cY$ is given and $F\in\cX$ is unknown. In our case, we
assume that $\tilde{G}\in\iota(\cZ)$, which might mean that only part
of the information, $\cP \tilde{G}\in\cZ$, of the \lq whole\rq\
right-hand side $\tilde{G}$ is given. For these reasons, we will deal
here with the inverse problem
\begin{equation}\label{Eq:IPdef}
    \cP \cT F= G, \quad G\in\cZ \text{ given, } F\in\cX \text{ unknown}\,.
\end{equation}
Unfortunately, we have the svd for the operator $\cT:\cX\to\cY$, but
not for the operator $\cP\cT:\cX\to\cZ$ . As we will see, the basic
principle of a Slepian approach is to obtain an svd for $\cP\cT$, which
is useful in cases where data are only obtainable from $\cZ$.

\begin{example}\label{Ex:upcontoperator}
We continue with the inverse problem from
Example \ref{Ex:upcontspaces}, the downward continuation problem
of \cite{PlattnerSimons2017}. The forward operator
$\cT\colon \LOp\to\LOs$ has the svd
\begin{equation}
\cT F= \sumlm \left(\frac{r_\mathrm{p}}{r_\mathrm{s}}\right)^l
\left\langle F, \frac{1}{r_\mathrm{p}}Y_{l,m}\left(\frac{\cdot}{r_\mathrm{p}}\right)\right\rangle_\LOplr\,
\frac{1}{r_\mathrm{s}} Y_{l,m}\left(\frac{\cdot}{r_\mathrm{s}}\right)\,,
\label{Eq:upcontsvd}
\end{equation}
for $F\in\LOp$, where $(Y_{l,m})_{l\in\N_0;\;m=-l,\dots,l}$ is the
commonly used orthonormal basis of real spherical harmonics in
$\rL^2(\Omega)$. Here, $R\subset\Os$ is the subdomain of data
availability, or of modeling interest for the potential.
\end{example}

\section{Slepian approach}\label{Sec:Slepian}
In analogy with \cite{PlattnerSimons2017}, we pursue the idea to
maximize
\begin{equation}
\cR(F)\coloneqq \frac{\left\|\cP \cT F\right\|^2_\cZ}{\|F\|_\cX^2}
\label{Eq:defR}
\end{equation}
among all $F\in\cX$ with $F\neq 0$. The individual terms can be
represented as follows ($\cP_{\ker \cT}$ is the orthogonal projection
onto the nullspace or kernel of $\cT$), for $F\in\cX$,
\begin{align}
\|F\|_\cX^2 &= \sum_n \left|\xip{F}{u_n}\right|^2 + \left\|\cP_{\ker \cT} F\right\|_\cX^2\,,\\
\cT F &= \sum_n \sigma_n \xip{F}{u_n} \,v_n\,,\label{Eq:svd2}\\
 \cP \cT F &= \sum_n \sigma_n \xip{F}{u_n} \,\cP v_n\,,\\
\|\cP\cT F\|_\cZ^2 &= \sum_{m,n} \sigma_m\overline{\sigma_n} \xip{F}{u_m}\overline{\xip{F}{u_n}}\zip{\cP v_m}{\cP v_n}.
\end{align}
Note that these formulae are also valid if $(v_n)_n$ is \emph{not}
orthonormal in $\cY$. It suffices that \eqref{Eq:svd}, which is the
same as \eqref{Eq:svd2}, is a finite sum or a (strongly)
convergent series.

\begin{example}
Let us consider again Example \ref{Ex:upcontoperator}. In this case,
the kernel is trivial, i.e.\ $\ker \cT =\{0\}$. Furthermore,
\begin{align}
\|F\|_\LOplr^2 &= \sumlm \left\langle F, \frac{1}{r_\mathrm{p}}Y_{l,m}\left(\frac{\cdot}{r_\mathrm{p}}\right)\right\rangle_\LOplr^2\,,\\
\|\cP\cT F\|_\LR^2 &= \sumlm \sum_{l'=0}^\infty \sum_{m'=-l'}^{l'} \left(\frac{r_\mathrm{p}}{r_\mathrm{s}}\right)^{l+l'}
\left\langle F, \frac{1}{r_\mathrm{p}}Y_{l,m}\left(\frac{\cdot}{r_\mathrm{p}}\right)\right\rangle_\LOplr\\
&\quad\times
\left\langle F, \frac{1}{r_\mathrm{p}}Y_{l',m'}\left(\frac{\cdot}{r_\mathrm{p}}\right)\right\rangle_\LOplr\\
&\quad \times \int_R \frac{1}{r_\mathrm{s}} Y_{l,m}\left(\frac{\eta}{r_\mathrm{s}}\right) \,
\frac{1}{r_\mathrm{s}} Y_{l',m'}\left(\frac{\eta}{r_\mathrm{s}}\right)\doe\,.
\end{align}
\end{example}

\begin{lemma}
The ratio $\cR$, which was defined in \eqref{Eq:defR}, satisfies
\begin{equation}
0\leq \cR (F)\leq \max_n|\sigma_n|
\end{equation}
for all $F\in\cX\setminus\{0\}$.
\end{lemma}
\begin{proof}
Since $\cP$ is a projection, its operator norm must satisfy
$\|\cP\|_{\cL(\cY,\cZ)}=1$. Moreover, the singular-value
decomposition \eqref{Eq:svd} yields that $\|\cT\|_{\cL(\cX,\cY)}
= \max_n|\sigma_n|$. Note that this maximum exists, since $\cT$ is
compact and, therefore, $(\sigma_n)_n$ must either be a finite
sequence or a sequence which converges to zero. Hence,
$\|\cP\cT\|_{\cL(\cX,\cZ)}\leq \max_n|\sigma_n|$, where
\begin{equation}
\|\cP\cT\|_{\cL(\cX,\cZ)}^2 =\sup_{F\in\cX\setminus\{0\}} \cR(F)\,.
\end{equation}
\end{proof}

\begin{lemma}
The operator $\cP\cT\colon\cX\to\cZ$ is compact.
\end{lemma}
\begin{proof}
$\cT$ is compact and $\cP$ is (as every projection) continuous. Hence, $\cP\cT$ is compact.
\end{proof}

As a consequence, $\cP\cT$ must have a singular-value decomposition
\begin{equation}
\cP\cT F= \sum_n \tau_n \xip{F}{g_n} \,h_n, \quad F\in\cX\,,
\label{Eq:svdSlep}
\end{equation}
where $(\tau_n)_n\subset\C$ is either a finite sequence or a sequence
converging to zero, $(g_n)_n$ is an orthonormal system in $\cX$, and
$(h_n)_n$ is an orthonormal system in $\cZ$. We will assume here that
the singular values $(\tau_n)_n$ are sorted in a way such that
$(|\tau_n|)_n$ is monotonically decreasing. The corresponding sequence
$(g_n)_n$ will be called a sequence of \textbf{Slepian basis
functions} with a localization of descending order. This is motivated
by the fact that
\begin{equation}
    \cR\left(g_n\right) = \frac{\left\|\cP\cT g_n\right\|_\cZ^2}{\left\|g_n\right\|_\cX^2} = |\tau_n|^2\,.
    \label{Eq:Rwithsvd}
\end{equation}

\section{Finite-dimensional case}\label{Sec:findim}
In numerical implementations, only finite basis systems can be
used. This usually means that the analysis is restricted to
bandlimited functions. Because of its practical relevance, we discuss
this particular case here separately. We set
\begin{align}
f&\coloneqq \left(\xip{F}{u_n}\right)_{n=1,\dots,N}\in\C^{N} \quad\text{(column vector)}\,, \label{Eq:findimf}\\
\Sigma&\coloneqq \begin{pmatrix}
\sigma_1 & 0 &\cdots & 0\\
0 & \ddots & \ddots &\vdots \\
\vdots & \ddots &\ddots & 0\\
0 &\cdots & 0 & \sigma_N
\end{pmatrix}\in\C^{N\times N}\, ,\label{Eq:findimA}\\
K&\coloneqq \left(\zip{\cP v_m}{\cP v_n}\right)_{m,n=1,\dots,N}\in\C^{N\times N}\, .\label{Eq:findimK}
\end{align}
Then the Parseval identity implies that
\begin{equation}
\cR(F)= \|f\|_{\C^{N}}^{-2} \cdot f^\mathrm{T} \Sigma K\Sigma^\ast \overline{f}\,.
\end{equation}
Here, $M^\ast$ represents the complex adjoint of a matrix
$M=(m_{i,j})_{i,j=1,\dots,N}$, i.e.\ $M^\ast \coloneqq
(\overline{m_{j,i}})_{i,j=1,\dots,N}$, $\overline{M}$ stands for the
complex conjugate
$\overline{M}\coloneqq(\overline{m_{i,j}})_{i,j=1,\dots,N}$, and
$M^\mathrm{T}\coloneqq(m_{j,i})_{i,j=1,\dots,N}$ is the transposed
matrix.

Since $f^\mathrm{T} \Sigma K\Sigma^\ast \overline{f}=\|\cP\cT F\|_\cZ^2$ is
real, $\Sigma$ is a diagonal matrix, and an inner product has a conjugate
symmetry, we can also write
\begin{equation}\label{Eq:findimR}
\cR(F)= \|f\|_{\C^{N}}^{-2} \cdot \overline{f^\mathrm{T} \Sigma K\Sigma^\ast
  \overline{f}} = \|f\|_{\C^{N}}^{-2} \cdot f^\ast \Sigma^\ast K^\mathrm{T}
\Sigma f\,,
\end{equation}
where
\begin{equation}\label{Eq:AKA}
    \Sigma^\ast K^\mathrm{T}\Sigma = \left(\overline{\sigma_m}\, \zip{\cP
      v_n}{\cP v_m} \,{\sigma_n}\right)_{m,n=1,\dots,N} \, .
\end{equation}
This result corresponds to the approach for the internal-field-only
case in \cite{PlattnerSimons2017}.

\begin{example}\label{Ex:downcontfindim}
We continue with Example \ref{Ex:upcontoperator}. In this case, the
index range is set to $l=0,\dots,L$, $m=-l,\dots,l$. Then the entries
of the diagonal matrix $\Sigma$ are given by
\begin{equation}
    \sigma_{l,m}= \left(\frac{r_\mathrm{p}}{r_\mathrm{s}}\right)^l\,.
\end{equation}
The unknown vector $f$ contains the Fourier coefficients
\begin{equation}
    f_{l,m}= \left\langle F,
    \frac{1}{r_\mathrm{p}}Y_{l,m}\left(\frac{\cdot}{r_\mathrm{p}}\right)\right\rangle_\LOplr\,.
\end{equation}
Furthermore, the matrix $K$ is given by its components
\begin{equation}
\int_R \frac{1}{r_\mathrm{s}} Y_{l,m}\left(\frac{\eta}{r_\mathrm{s}}\right) \,
\frac{1}{r_\mathrm{s}} Y_{l',m'}\left(\frac{\eta}{r_\mathrm{s}}\right)\doe\,.
\end{equation}
The task is, therefore, to find the eigenvectors $f$ of the matrix
\begin{equation}
\Sigma^\ast
K^\mathrm{T}\Sigma=\left[\left(\frac{r_\mathrm{p}}{r_\mathrm{s}}\right)^l\,\int_R
\frac{1}{r_\mathrm{s}} Y_{l,m}\left(\frac{\eta}{r_\mathrm{s}}\right)
\, \frac{1}{r_\mathrm{s}}
Y_{l',m'}\left(\frac{\eta}{r_\mathrm{s}}\right)\doe\,
\left(\frac{r_\mathrm{p}}{r_\mathrm{s}}\right)^{l'}
\right]_{\genfrac{}{}{0pt}{}{l=0,\dots,L;\;
    m=-l,\dots,l}{l'=0,\dots,L;\; m'=-l',\dots,l'}}\,.
\end{equation}
\end{example}

\begin{rem}
For some problems, vectorial (e.g.\ gradients of potential fields
\cite{PlattnerSimons2017}) or tensorial basis functions come into
play, and the $\rL^2$ inner products involve Euclidean dot products of
the kind
\begin{equation}
\zip{f}{g}=\int_R f(\xi)\cdot g(\xi) \dox,\quad f,g\in\rL^2\left(R,\R^3\right)\,.
\end{equation}
In this case, one can make use of this Euclidean product to reduce the
numerical expense or the instability of the eigenvalue problem at
hand. For example, the vector spherical harmonics (for which we use
here the notation in \cite{OxfordBuch}) $y_{l,m}^{(i)}$ can be
subdivided into vector fields which are normal to the sphere ($i=1$)
and fields that are tangential to the sphere ($i=2$ and $i=3$). For
this reason,
\begin{equation}
    y_{l,m}^{(1)}(\xi) \cdot y_{n,j}^{(i)}(\xi) = 0
\end{equation}
holds pointwise (i.e.\ for all $\xi\in\Omega$) and all $i\in\{2;3\}$,
independently of the degrees $l,n$ and orders $m,j$. Within the
tangential vector fields, such a pointwise, i.e.\ Euclidean,
orthogonality is only obtained for identical degree-order pairs, i.e.
\begin{equation}
y_{l,m}^{(2)}(\xi) \cdot y_{l,m}^{(3)}(\xi) = 0
\end{equation}
for all $\xi\in\Omega$ and all degrees $l$ and orders $m$. 

In \cite{JahnBokor2014}, different linear combinations of
complex tangential vector spherical harmonics are constructed to
obtain alternative basis functions, which we call here
$\tilde{y}_{l,m}^{(i)}$, $i=1,2,3$, such that\footnote{Note that
  \lq$\cdot$\rq\ is here the complex dot product, i.e.\ $w\cdot
  z\coloneqq \sum_{j=1}^3 w_j \overline{z_j}$ for $w,z\in\C^3$.}
\begin{equation}
    \tilde{y}_{l,m}^{(2)}(\xi) \cdot \tilde{y}_{n,j}^{(3)}(\xi) = 0
\end{equation}
for all $\xi\in\Omega$, all degrees $l,n$, and all orders $m,j$. This pointwise orthogonality can be exploited, because we have
\begin{equation}
\int_R \tilde{y}_{l,m}^{\left(i_1\right)}(\xi) \cdot \tilde{y}_{n,j}^{\left(i_2\right)}(\xi) \dox =0
\end{equation}
whenever $i_1\neq i_2$. As a consequence, the matrix $K$ can be rearranged into a block matrix
\begin{equation}
\begin{pmatrix}
\text{type } i=1 & 0 & 0 \\
0 & \text{type } i=2 & 0 \\
0 & 0 & \text{type } i=3
\end{pmatrix}
\end{equation}
such that the algebraic eigenvalue problems can be solved separately
for each type $i$. This has not only the advantage that the matrices
of the eigenvalue problem become smaller (which yields the expectation
of a faster and more stable computation of the eigenvectors), it also
leads to Slepian functions which are separated by type. This means
that the components of the field associated to different types can be
independently analyzed by means of Slepian functions. 

However, one has to be aware of the fact that the type $i=2$ of the
$y_{l,m}^{(i)}$ (which is a surface gradient field and is, therefore,
surface-curl-free) is not the same as the type $i=2$ of the
$\tilde{y}_{l,m}^{(i)}$. The reason is that each tangential
$\tilde{y}_{l,m}^{(i)}$, $i\in\{2;3\}$, is a linear combination of
complex versions of $y_{l,m}^{(2)}$ \emph{and} $y_{l,m}^{(3)}$. In
particular, $\tilde{y}_{l,m}^{(2)}$ is not surface-curl-free anymore,
and $\tilde{y}_{l,m}^{(3)}$ is not surface-divergence-free anymore ---
properties which the non-tilde versions originally possessed.

For tensor spherical harmonics, there are 9 different types of basis
functions, where again some types are orthogonal to each other in the
Euclidean sense. Also here, it is possible to define a new basis
system such that the Slepian eigenvalue problem can be transformed
into 9 independent eigenvalue problems, as shown in
\cite{SeibertDiss2017}.
\end{rem}

\section{Coupled problems}\label{Sec:coup}
In some applications, we may have data that originate from different
causes or sources, and we may be interested in separating them (e.g.\
internally and externally generated planetary magnetic
fields \cite{PlattnerSimons2017}). We will show here that such a
scenario can easily be integrated into our general setting.

We now have two operators $\cT_1\colon\cX_1\to\cY$ and $\cT_2\colon\cX_2\to\cY$  with svds
\begin{align}
\cT_1 F_1 &= \sum_n \sigma_n^{(1)} \xoip{F_1}{u_n^{(1)}} \, v_n^{(1)}, \quad F_1\in\cX_1\,,\\
\cT_2 F_2 &= \sum_n \sigma_n^{(2)} \xtip{F_2}{u_n^{(2)}} \, v_n^{(2)}, \quad F_2\in\cX_2\, .
\end{align}
The notation for the Hilbert spaces and the orthonormal systems is
analogous to the previous case. Note that $\cT_1$ and $\cT_2$ both map
into $\cY$, but that they may use different orthonormal systems
$(v_n^{(1)})_n\subset\cY$ and $(v_n^{(2)})_n\subset\cY$. As in the
single-operator case, we have only one Hilbert space $\cZ$, one
projection $\cP\colon\cY\to\cZ$, and one embedding
$\iota\colon\cZ\hookrightarrow\cY$.

The inverse problem is now to find $F_1\in\cX_1$ and $F_2\in\cX_2$
such that, for a given $G\in\cZ$,
\begin{equation}
\cP\cT_1F_1+\cP\cT_2F_2 = G
.
\end{equation}

\begin{example}
In \cite{PlattnerSimons2017}, it is assumed that a potential field is
given which is a superposition of potentials from an internal and an
external source, where the sources could be of a magnetic or a
gravitational nature. More precisely, the case of gradients of the
potential is considered. For reasons of brevity of the formulae, we
will consider here the scalar potential situation. The inner potential
corresponds to Example \ref{Ex:upcontoperator} and its source is
assumed to be located inside the planet (i.e.\ in the interior of
$\Op$). The external potential originates from a radius of at least
$r_\mathrm{e}$, where $r_\mathrm{e}>r_\mathrm{s}$. This leads to the
operators
\begin{subequations}\label{Eq:couple}
\begin{align}
\cT_1 F_1&= \sumlm \left(\frac{r_\mathrm{p}}{r_\mathrm{s}}\right)^l
\left\langle F_1, \frac{1}{r_\mathrm{p}}Y_{l,m}\left(\frac{\cdot}{r_\mathrm{p}}\right)\right\rangle_\LOplr\,
\frac{1}{r_\mathrm{s}} Y_{l,m}\left(\frac{\cdot}{r_\mathrm{s}}\right)\,,\label{Eq:couple_1}\\
\cT_2 F_2&= \sumlm \left(\frac{r_\mathrm{s}}{r_\mathrm{e}}\right)^{l+1}
\left\langle F_2, \frac{1}{r_\mathrm{e}}Y_{l,m}\left(\frac{\cdot}{r_\mathrm{e}}\right)\right\rangle_\LOelr\,
\frac{1}{r_\mathrm{s}} Y_{l,m}\left(\frac{\cdot}{r_\mathrm{s}}\right)\,,\label{Eq:couple_2}
\end{align}
\end{subequations}
where \eqref{Eq:couple_1} represents the inner field
and \eqref{Eq:couple_2} stands for the external field. The singular
values of $\cT_1$ and $\cT_2$ both exponentially converge to $0$,
which means that both operators are compact. However, depending on the
values of $r_\mathrm{p}$, $r_\mathrm{s}$, and $r_\mathrm{e}$, these
two sequences need not tend to zero equally fast. This means that the
associated ill-posednesses need not be equally severe. As a
consequence, it can be reasonable to truncate the two series
in \eqref{Eq:couple} at different degrees. The consequently different
sizes of the orthonormal systems combined with the different
instabilities (and, maybe also coupled with different noise scenarios)
yield a situation which can be expected to be particularly challenging
regarding the necessary regularization.
\end{example}

Let us return to the general setting. For the considered problem, we
construct the Hilbert space $\cX\coloneqq \cX_1\otimes \cX_2$ as the
Cartesian product of the individual spaces, and equip it with the
inner product, for $x_1,x_1'\in\cX_1$, $x_2,x_2'\in\cX_2$,
\begin{equation}
\xip{(x_1,x_2)}{(x_1',x_2')} \coloneqq \xoip{x_1}{x_1'} + \xtip{x_2}{x_2'}\,,
\end{equation}

Moreover, we define the operator $\cS\colon \cX\to\cY$ by
\begin{equation}
\cS \left(F_1,F_2\right) \coloneqq \cT_1F_1+\cT_2F_2, \quad F_1\in\cX_1,\, F_2\in\cX_2\, .
\end{equation}

Furthermore, we set
\begin{align}
u_{2n} &\coloneqq \left(u_n^{(1)},0\right), & u_{2n+1}&\coloneqq \left(0,u_n^{(2)}\right)\,, \\
v_{2n} &\coloneqq v_n^{(1)}, & v_{2n+1}&\coloneqq v_n^{(2)}\,,\\
\sigma_{2n}&\coloneqq \sigma_n^{(1)}, & \sigma_{2n+1}&\coloneqq\sigma_n^{(2)}\,.
\end{align}
This arrangement of the two systems into one system certainly does not
necessarily have to be done in this order. In particular, in the
finite-dimensional case, where we only have
$(u_n^{(1)})_{n=1,\dots,N_1}$ and $(u_n^{(2)})_{n=1,\dots,N_2}$, we
could equivalently set
\begin{equation}
\left(u_1,\dots,u_{N_1+N_2}\right)\coloneqq\left(u_1^{(1)},\dots,u_{N_1}^{(1)},u_{1}^{(2)},\dots,u_{N_2}^{(2)}\right)\,.
\end{equation}
We now have
\begin{equation}
\cS F = \sum_n \sigma_n \xip{F}{u_n} \, v_n,\quad F\in\cX\,,
\end{equation}
where $(u_n)_n$ is an orthonormal system in $\cX$ \emph{but} $(v_n)_n$
is, in general, \emph{not} an orthonormal system in $\cY$. In
Section \ref{Sec:Slepian} we remarked that there is no requirement
that $(v_n)_n$ be orthonormal, hence we can proceed now like in
the \lq non-coupled\rq\ case. However, in the (theoretical) case where
infinite systems are involved, the particular arrangement of the two
systems into one system could be of importance in the sense of the
Riemann series theorem (see e.g.\ \cite[p.\ 68]{Bromwich1908})

In the finite-dimensional case, the Slepian matrix,
\begin{equation}
    \Sigma^\ast K^\mathrm{T}\Sigma=\left(\overline{\sigma_m}\, \zip{\cP v_n}{\cP v_m} \,{\sigma_n}\right)_{m,n=1,\dots,N}
\end{equation}
corresponds to the matrix of the eigenvalue problem for the
mixed-source case in \cite{PlattnerSimons2017}.

Couplings of more than two sources can be handled analogously.

\section{Solving the inverse problem}\label{Sec:solv_ip}
\subsection{The Slepian functions and the svd}
The svd for the operator $\cP\cT\colon\cX\to\cZ$ in
\eqref{Eq:svdSlep}, which we know exists, allows us to use a truncated
singular-value decomposition
\begin{equation}
F_J=\sum_{\genfrac{}{}{0pt}{}{k=1}{\tau_k\neq 0}}^J \tau_k^{-1} \zip{G}{h_k} \, g_k\label{Eq:TSVDSlep}
\end{equation}
as an approximate solution of the inverse problem $\cP\cT F=G$,
$G\in\cZ$.

To find $g_k$, $h_k$, and $\tau_k$, in the finite-dimensional setting
of Section \ref{Sec:findim}, the procedure requires us to:
\begin{itemize}
  \item set up the matrix $\Sigma^\ast K^\mathrm{T} \Sigma$ as in \eqref{Eq:AKA}.
  \item determine an orthonormal system of eigenvectors\footnote{Since
  $\Sigma^\ast K^\mathrm{T}\Sigma$ is self-adjoint, such an orthonormal basis
  must exist and all eigenvalues are real. Moreover, since it is a
  Gramian matrix, all eigenvalues must be non-negative.}
  $f^{(1)},\dots,f^{(N)}$ and its associated eigenvalues
  $\varrho_1,\dots,\varrho_{N}\in\R_0^+$. 
  \item sort the eigenvalues (and the associated eigenvectors) such that
  $\varrho_1\geq\varrho_2\geq\dots\geq\varrho_{N}$.
  \item construct
  \begin{equation}
    g_k=\sum_{n=1}^N f^{(k)}_n u_n \in\cX, \quad k=1,\dots, N\,
    ,
  \end{equation}
  with the Parseval identity yielding
  \begin{equation}
  \xip{g_k}{g_l} = \sum_{n=1}^N\ f^{(k)}_n \overline{f^{(l)}_n} = \left\langle f^{(k)},f^{(l)} \right\rangle_{\C^{N}} = \delta_{kl}\,. 
  \end{equation}
\end{itemize}
We now use the $g_k$ as basis functions to expand the solution $F$,
that is, we determine coefficients $\gamma_k$ such that
$F=\sum_{k=1}^N \gamma_k g_k$ solves $\cP\cT F=G$. In keeping with the
common philosophy of Slepian functions, we may truncate our expansions
by taking only these Slepian functions $g_k$ for which
$\varrho_k\geq\tilde{\varrho}$ for a chosen threshold
$\tilde{\varrho}$.

In our case, we determine the $\gamma_k$ from the svd of $\cP\cT$,
proceeding as follows. From \eqref{Eq:svd}, we know
\begin{equation}
\cT g_k = \sum_{n=1}^N \sigma_n f^{(k)}_n v_n, \quad k=1,\dots,N\,,
\end{equation}
and therefore also
\begin{equation}
\cP\cT g_k = \sum_{n=1}^N \sigma_n f^{(k)}_n \cP v_n, \quad k=1,\dots,N
.
\end{equation}
Furthermore, with \eqref{Eq:AKA}, an interchanging of $m$ and $n$, and
the fact that the $f^{(k)}$ are orthonormal eigenvectors of $\Sigma^\ast
K^\mathrm{T}\Sigma$, we get 
\begin{align}
\zip{\cP\cT g_k}{\cP\cT g_l} 
&= \sum_{m,n=1}^N \sigma_m\overline{\sigma_n} f_m^{(k)}\overline{f_n^{(l)}} \zip{\cP v_m}{\cP v_n}\\
&= {f^{(l)}}^\ast \Sigma^\ast K^\mathrm{T}\Sigma \,{f^{(k)}}\\
&= \varrho_k \overline{\left\langle f^{(l)}, f^{(k)}\right\rangle_{\C^{N}}}\\
&=\varrho_k \delta_{k,l}\,.
\end{align}
We set
\begin{align}
h_k&\coloneqq \varrho_k^{-1/2} \cP\cT g_k, \quad k=1,\dots,N\,,\quad \text{if }\varrho_k\neq 0\,,\\
\tau_k&\coloneqq \varrho_k^{1/2}, \quad k=1,\dots,N,\quad \text{if }\varrho_k\neq 0\,.
\end{align}
If $\varrho_k=0$, then $\cP\cT g_k=0$ such that we set
$\tilde{N}\coloneqq\max\{k\,|\,\varrho_k\neq0\}$ and consider
$(u_k)_{k=1,\dots,\tilde{N}}$ as an orthonormal basis of
$(\ker\cP\cT)^{\bot_\cX}$, the $\xip{\cdot}{\cdot}$-orthogonal complement of the nullspace
of $\cP\cT$.

Then,
\begin{align}
\cP\cT F&= \cP\cT \sum_{k=1}^{\tilde{N}} \xip{F}{g_k} g_k\\
&= \sum_{k=1}^{\tilde{N}} \xip{F}{g_k} \tau_k\,h_k\,,\quad F\in\cX.
\end{align}
This is the required svd of $\cP\cT\colon\cX\to\cZ$, see \eqref{Eq:svdSlep}.

The determination of the truncation parameter $J$
in \eqref{Eq:TSVDSlep} can be accomplished with any one of the known
parameter choice methods for the regularization of inverse problems (see e.g.\ \cite{Gutting+2017,Xu98} and the references therein).
Furthermore, $(h_k)_{k=1,\cdots,\tilde{N}}$ is an orthonormal system
in $\cZ$. Moreover, \eqref{Eq:isometry} implies that
\begin{equation}\label{Eq:hkonsisom}
\delta_{kl}=\zip{h_k}{h_l} = \yip{\iota(h_k)}{\iota(h_l)}
\end{equation}
such that $(\iota(h_k))_{k=1,\dots,\tilde{N}}$ is also orthonormal in
$\cY$. In our example of function spaces, this means that the $\iota(h_k)$
are spacelimited functions, which are orthogonal in $\LOs$.

\begin{example}
We continue with Example \ref{Ex:downcontfindim}, reverting to the
degree and order indices $l=0,\dots,L$, $m=-l,\dots,l$. After having
obtained the eigenvectors $f^{(k)}$ and eigenvalues $\varrho_k$ for
the matrix $\Sigma^\ast K^\mathrm{T}\Sigma$, we can calculate the
following functions:
\begin{align}
g_k(\xi)&= \sumlmL f_{l,m}^{(k)} \,\frac{1}{r_\mathrm{p}}\, Y_{l,m}\left(\frac{\xi}{r_\mathrm{p}}\right) , \quad \xi\in\Op\,,\\
h_k(\zeta)&= \varrho_k^{-1/2} \sumlmL f_{l,m}^{(k)} \left(\frac{r_\mathrm{p}}{r_\mathrm{s}}\right)^l \frac{1}{r_\mathrm{s}}\,Y_{l,m}\left(\frac{\zeta}{r_\mathrm{s}}\right), \quad \zeta\in R\,,\\
\left[\iota\left(h_k\right)\right](\eta) &= \begin{cases}
\varrho_k^{-1/2} \sumlmL f_{l,m}^{(k)} \left(\frac{r_\mathrm{p}}{r_\mathrm{s}}\right)^l \frac{1}{r_\mathrm{s}}\,Y_{l,m}\left(\frac{\eta}{r_\mathrm{s}}\right),& \eta\in R\\
0,& \eta\in\Os\setminus R
\end{cases}
, \quad \eta\in\Os\,.
\end{align}
We obtain then the following orthogonalities
\begin{align}
\int_\Op g_k (\xi)\, g_l(\xi)\dox &= \delta_{kl},\\
\int_R h_k(\zeta)\,h_l(\zeta)\doz
&= \int_\Os \left[\iota\left(h_k\right)\right](\eta) \,\left[\iota\left(h_l\right)\right](\eta)\doe
= \delta_{kl}\,. 
\end{align}
\end{example}
The example and the considerations above show one of the advantages of
Slepian functions applied to inverse problems with regional data. We
are able to obtain a singular-value decomposition for the projected
operator $\cP\cT\colon\cX\to\cZ$, that is, for the case where only
regional data are available. We have orthonormal function systems
$(g_k)_k$ in~$\cX$ and $(h_k)_k$ in~$\cZ$ which can be calculated
explicitly.

\subsection{Construction of a scaling function as a filter}

Moreover, alternative methods like wavelet-based multiscale methods
are applicable, where we introduce a filter $ \varphi_J$, such that
\begin{equation}\label{Eq:scalfct}
\tilde{F}_J=\sum_{k=1}^{\tilde{N}} \varphi_J(k) \tau_k^{-1} \zip{G}{h_k} \, g_k\,.
\end{equation}
In the case of functions, this could be
\begin{equation}\label{Eq:scalfctfct}
\tilde{F}_J(x) = \zip{G}{\Phi_J(x,\cdot)}\, ,
\end{equation}
where the scaling function $\Phi_J$ is given by
\begin{equation}\label{Eq:scalfctformula}
\Phi_J(x,z)=\sum_{k=1}^{\tilde{N}} \varphi_J(k) \,\tau_k^{-1} \,\overline{g_k(x)}\,h_k(z)\,.
\end{equation}
We will further elaborate this in Section~\ref{Sec:Fredholm}.

\subsection{Infinite-dimensional case}

Putting numerical considerations aside for a moment, we can observe
that the considerations here are not restricted to the
finite-dimensional case. With the (initially unknown but definitely
existing) singular-value decomposition \eqref{Eq:svdSlep} and with
\eqref{Eq:Rwithsvd}, we could also proceed with an infinite (e.g.\ non-bandlimited) setting. We would get a (possibly infinite, but
countable) system of non-negative values
$(\tau_k)_{k\in\kappa}=(\varrho_k^{1/2})_{k\in\kappa}$, where
$\kappa\subset\N$ stands here for the index range which counts all
such singular values. Due to the nature of an svd, the $g_k$,
$k\in\kappa$, would represent an orthonormal system in $\cX$. More
precisely, we would have an orthonormal basis of
$(\ker\cP\cT)^{\bot_\cX}$. Then
\begin{equation}
h_k= \tau_k^{-1} \cP\cT g_k ,\quad k\in\kappa\,,
\end{equation}
is an orthonormal system such that $\mathop{\mathrm{im}}
\cP\cT\subset\overline{\mathop{\mathrm{span}}\{h_k\,|\,k\in\kappa\}}^{\|\cdot\|_\cZ}$
(the closure of the span of $h_k$, i.e.  every element in the image of
$\cP\cT$ can be expanded into the basis $h_k$, possibly with an
infinite number of summands). Furthermore, we would also get that
$(\iota(h_k))_{k\in\kappa}$ is orthonormal in $\cY$ due to $\iota$
being an isometry.

\section[Scaling functions, wavelets etc.]{Scaling functions, wavelets, reproducing kernels and Fredholm integral operators}\label{Sec:Fredholm}
In this section, we assume that the Hilbert spaces $\cX$, $\cY$, and
$\cZ$ are spaces of functions with domains $X$, $Y$, and $Z$,
respectively. The example of downward continuation which has been
discussed throughout this paper fits this assumption. 

For an $F\in\cX$, a $G\in\cZ$, we use the svd of the problem $\cP\cT F=G$,
\begin{align}
(\cP\cT F)(z) &= \sum_{k\in\kappa} \tau_k \xip{F}{g_k} \, h_k(z)\\
&= \xip{F(\cdot)}{\sum_{k\in\kappa} \tau_k \,g_k(\cdot)\,\overline{h_k(z)}}\\
&=\xip{F(\cdot)}{D^\uparrow(z,\cdot)},
\end{align}
to define now the following functions
\begin{align}
D^\uparrow (z,x)&\coloneqq \sum_{k\in\kappa} \tau_k\,\overline{h_k(z)}\,g_k(x), \quad z\in Z,\, x\in X\,,\\
D^\downarrow (x,z)&\coloneqq \sum_{k\in\kappa} \tau_k^{-1}\,\overline{g_k(x)}\,h_k(z), \quad x\in X,\, z\in Z\,,
\end{align}
assuming appropriate convergence\footnote{We need that, for each fixed
$z\in Z$, the series corresponding to $D^\uparrow(z,\cdot)$ converges
strongly in the sense of $\|\cdot\|_\cX$. Analogously, for each fixed
$x\in\cX$, the series corresponding to $D^\downarrow(x,\cdot)$ must be
strongly convergent in the sense of $\|\cdot\|_\cZ$.} in the case of
an infinite number of summands.  Similarly,
\begin{align}
\left([\cP\cT]^{+}G\right) (x) &= \sum_{k\in\kappa} \tau_k^{-1} \zip{G}{h_k} \, g_k(x)\\
&= \zip{G(\cdot)}{\sum_{k\in\kappa} \tau_k^{-1} \, h_k(\cdot)\,\overline{g_k(x)}}\\
&=\zip{G(\cdot)}{D^\downarrow(x,\cdot)},\quad x\in X\,,
\end{align}
where $(\cP\cT)^{+}=(\cP\cT|_{(\ker(\cP\cT))^{\bot_\cX}})^{-1}$ is the Moore-Penrose inverse of $\cP\cT$.

The kernel $D^\downarrow$ probably will not exist in the
infinite-dimensional case, because $(\tau_k^{-1})$ diverges to
$+\infty$. This represents the ill-posedness of the problem, because
$(\cP\cT)^{+}G$ cannot so easily be computed. For this reason, a
regularization is needed. 

This can be done in manifold ways, where a truncation of the series,
which would be the classical Slepian function approach discussed above
in~\eqref{Eq:TSVDSlep}, is one out of these possibilities. The more
general Ansatz corresponds to the scaling function approach described
above in~\eqref{Eq:scalfctformula}, where we replace $D^\downarrow$
by the kernel
\begin{equation}
\Phi_J(x,z)=\sum_{k\in\kappa} \varphi_J(k) \,\tau_k^{-1} \,\overline{g_k(x)}\,h_k(z)\,, \quad x\in X,\, z\in Z\,. 
\end{equation}
By choosing a sequence $(\varphi_J(k))_k$, which tends to zero \lq
sufficiently\rq\ fast, we can control the rising inverse singular
values $\tau_k^{-1}$ and obtain a stable solution. Such wavelet-based
regularization methods have already been discussed for such general
Hilbert space settings in \cite{FreGloLit1999,MichelHab}. We will show
here the most important properties of such a multiscale regularization
for the considered Slepian-function approach.
\begin{theorem}\label{Th:regul}
Let the assumptions from above hold true. Moreover, let the family of
functions $\varphi_J\colon\R^+_0\to\R_0^+$, $J\in\N_0$, satisfy the
following conditions\footnote{If $\kappa$ is a finite set, then
conditions (\ref{it:Th_regul_l2conv}) and (\ref{it:Th_regul_condstab})
are trivially satisfied.}:
\begin{enumerate}[(i)]
  \item\label{it:Th_regul_l2conv} for all $J\in\N_0$ and all $x\in X$, the following series converges pointwise:
  \begin{equation}
    \sum_{k\in\kappa} \left|\varphi_J(k)\,\overline{g_k(x)}\,\tau_k^{-1}\right|^2<+\infty\,,
    \label{Eq:regul1}
  \end{equation}
  \item\label{it:Th_regul_condstab}for all $J\in\N_0$,
      \begin{equation}
      \sup_{k\in\kappa} \left(\varphi_J(k)\tau_k^{-1}\right) <+\infty\,,
      \label{Eq:Th_regul_condstab}
      \end{equation}
  \item for all $k\in\kappa$, 
  \begin{equation}
  \lim_{J\to\infty} \varphi_J(k) = 1\,,
  \label{Eq:regul2}
  \end{equation}
  \item for all $J\in\N_0$ and all $k\in\kappa$,
  \begin{equation}
  0\leq \varphi_J(k) \leq 1\,,
  \label{Eq:regul3}
  \end{equation}
\end{enumerate}
Furthermore, the sequence of functions $\Phi_J\ast G\in\cX$, $J\in\N_0$, is defined by
\begin{equation}\label{Eq:defconv}
    \left(\Phi_J\ast G\right)(x) \coloneqq \zip{G(\cdot)}{\sum_{k\in\kappa} \varphi(k)\,\tau_k^{-1} \,\overline{g_k(x)}\,h_k(\cdot)},\quad x\in X,\, G\in\cZ\,.
\end{equation}
Then
\begin{equation}
\lim_{J\to\infty}\left\|[\cP\cT]^{+} G - \Phi_J\ast G\right\|_\cX = 0
\end{equation}
for all $G\in\mathop{\mathrm{im}}(\cP\cT)$. Moreover, each mapping
\begin{align}
\cZ&\to \cX\\
G& \mapsto \Phi_J\ast G\,,
\end{align}
$J\in\N_0$, is continuous.
\end{theorem}
Before we prove this theorem, let us state what it means for the
inverse problem. The sequence $(\Phi_J\ast G)_J$ converges strongly
(in the $\|\cdot\|_\cX$-sense) to the solution $F\in(\ker
(\cP\cT))^{\bot_\cX}$ of the inverse problem $\cP\cT F=G$, provided
that a solution exists (i.e.\
$G\in \mathop{\mathrm{im}}(\cP\cT)$). Hence, we can construct
approximate solutions which are arbitrarily close to the exact
solution. However, in contrast to the exact solution $F$, which
discontinuously depends on $G$ in the infinite-dimensional case
(remember that $\cP\cT$ is compact), the approximations are stable,
that is they continuously depend on the data $G$. This also yields the
expectation of numerically stable approximate inversions in the
finite-dimensional case.

Let us now prove the theorem.
\begin{proof}
From the condition in \eqref{Eq:regul1}, we obtain that the series
\begin{equation}
\sum_{k\in\kappa} \varphi_J(k)\,\tau_k^{-1} \,\overline{g_k(x)}\,h_k(\cdot)\,,
\end{equation}
with arbitrary but fixed $J\in\N$ and $x\in X$, converges strongly in
$\cZ$. Hence, we are allowed to interchange the inner product with the
series in \eqref{Eq:defconv} and get
\begin{equation}
\left(\Phi_J\ast G\right)(x) = \sum_{k\in\kappa} \varphi_J(k)\, \tau_k^{-1}\, \zip{G}{h_k}\,g_k(x)
\end{equation}
for all $J\in\N_0$ and all $x\in X$. Furthermore, the solvability of the
inverse problem $\cP\cT F=G$ yields a unique (minimum-norm) solution
$F\in(\ker (\cP\cT))^{\bot_\cX}$, which is given by
\begin{equation}
F= \sum_{k\in\kappa} \tau_k^{-1}\, \zip{G}{h_k}\,g_k
\end{equation}
in the sense of $\|\cdot\|_\cX$. Hence, the well-known Picard condition
\begin{equation}
\sum_{k\in\kappa} \left| \tau_k^{-1}\, \zip{G}{h_k}\right|^2 <+\infty
\end{equation}
must hold. This Picard condition in combination with \eqref{Eq:regul3}
implies that the series
\begin{equation}
\left\| F-\Phi_J\ast G\right\|_\cX^2 =\sum_{k\in\kappa} \left|\left(1-\varphi_J(k)\right)\tau_k^{-1}\, \zip{G}{h_k}\right|^2
\end{equation}
uniformly converges with respect to all $J\in\N$. Hence,
\begin{align}
\lim_{J\to\infty}\left\| F-\Phi_J\ast G\right\|_\cX^2 &=\lim_{J\to\infty}\sum_{k\in\kappa} \left|\left(1-\varphi_J(k)\right)\tau_k^{-1}\, \zip{G}{h_k}\right|^2\\
&=\sum_{k\in\kappa} \left|\left(1-\lim_{J\to\infty}\varphi_J(k)\right)\tau_k^{-1}\, \zip{G}{h_k}\right|^2\\
&= 0
\end{align}
due to \eqref{Eq:regul2}.

For proving the stability of the approximations $\Phi_J\ast G$, we have a look at
\begin{align}
\left\|\Phi_J\ast G\right\|_\cX^2 & = \sum_{k\in\kappa} \left|\varphi_J(k)\,\tau_k^{-1}\,\zip{G}{h_k}\right|^2\notag\\
&\leq\sup_{k\in\kappa} \left(\varphi_J(k)\,\tau_k^{-1}\right)^2 \,\sum_{k\in\kappa} \left|\zip{G}{h_k}\right|^2\notag\\
&=\sup_{k\in\kappa} \left(\varphi_J(k)\,\tau_k^{-1}\right)^2 \,\|G\|_\cZ^2,\quad G\in\cZ\,. \label{Eq:provestab}
\end{align}
According to \eqref{Eq:Th_regul_condstab}, the supremum in \eqref{Eq:provestab} is finite.
This proves the continuity of the mapping.
\end{proof}

Examples for the choice of $\varphi_J$ can be constructed out of
generators of scaling functions as they are known, for instance, from
the theory of spherical wavelets (see e.g.\ \cite[Sections 11.3 and
11.4]{OxfordBuch}, \cite[Example 2.3.7]{MichelHab}, and \cite[Example
7.20]{Michelbook2013}). However, the critical part is represented by
conditions (\ref{it:Th_regul_l2conv}) and
(\ref{it:Th_regul_condstab}). They can be trivially satisfied by
taking generators of bandlimited scaling functions, that is functions
$\varphi_J$ with compact support $\supp \varphi_J$ for each
$J\in\N_0$. In the non-bandlimited case, where the support is
unbounded for an infinite number of scales $J$, the particular
properties of the computed Slepian functions $g_k$ and the rate of
divergence of the inverse singular values $(\tau_k^{-1})$ have to be
taken into account. On the one hand, this yields an interesting
challenge for future research, because these requirements implicitly
also include the geometry of the region $R$ (as well as the degree of
the ill-posedness of the original inverse problem $\cT F=G$) into the
conditions on $\varphi_J$. On the other hand, in practice, one either
always has to restrict the calculations to finite dimensional spaces,
that is, to the bandlimited case, or $\kappa$ is a finite set to begin
with.

Note also that, in the particular case of $\rL^2$-inner products in
$\cX$ and $\cZ$, we can, indeed, write the inverse problem as a
Fredholm integral equation of the first kind
\begin{equation}\label{Eq:Fredie1}
(\cP\cT F)(z) =\int_X F(x)\,\overline{D^\uparrow(z,x)}\dx, \quad z\in Z\,.
\end{equation}

Let us discuss now a special case: $\cT=\cI$ (identity) and $\cX=\cY$,
i.e.\ we \lq simply\rq\ want to interpolate/approximate a function. In
this case, $u_n=v_n$ for all $n$ and $\sigma_n=1$ for all $n$. The
singular-value decomposition of $\cT$ would be representable as
\begin{equation}
\cT F= \sum_n \xip{F}{u_n}\, u_n, \quad F\in\cX\,.
\end{equation}
The task is still to find a new singular-value decomposition for the
projected equation, but this time it is only the projection itself
which needs the svd. We, therefore, look for a representation of the
form
\begin{equation}
\cP F= \sum_k \tau_k \xip{F}{g_k} \, h_k
\end{equation}
which originates in the same way from the eigenvalue- or
singular-value-problem discussed above, where now
$h_k=\varrho_k^{-1/2}\cP g_k=\tau_k^{-1}\cP g_k$. If $\cP$ is the
restriction operator $\cP\colon F\mapsto \left.F\right|_Z$, then
\begin{align}
D^\uparrow(z,x) &= \sum_{k\in\kappa} \overline{g_k(z)}\,g_k(x), \quad z\in Z,\, x\in X\,,\\
D^\downarrow(x,z)&= \sum_{k\in\kappa} \tau_k^{-2}\,\overline{g_k(x)}\,g_k(z), \quad x\in X,\, z\in Z\,.
\end{align}
In other words, using again $\rL^2$-inner products, we see that
\begin{equation}
\int_X F(x) \overline{D^\uparrow (z,x)}\dx = (\cP F)(z)=F(z), \quad z\in Z\,,
\end{equation}
\emph{reproduces} $F$ on the subset $Z\subset Y=X$. In particular,
\begin{equation}
\int_X g_k(x) \overline{D^\uparrow (z,x)}\dx = (\cP g_k)(z)=g_k(z), \quad z\in Z\,.
\end{equation}
Vice versa,
\begin{equation}
\int_Z F(z) \overline{D^\downarrow (x,z)}\dz = \left(\cP^{+} F\right)(x), \quad x\in X\,,
\end{equation}
\emph{reconstructs} $F$ on the whole set $X$ from knowledge of $F$ on
the subset $Z$. The latter sounds confusing at the first sight. How
could the continuation of $F$ to a larger set be unique? Indeed, there
is a catch: the series of $D^\downarrow$ must converge. This can be
satisfied in two cases:   
\begin{itemize}
  \item either $\kappa$ is finite: then the function spaces under
  investigation have finite dimensions and the functions in it are,
  indeed, uniquely determined by their values on a subset (like it is
  e.g.\ the case for polynomials up to a fixed degree), 
\item or   $\kappa$ is infinite but the series converges nevertheless: then
  this implies certain regularity conditions on the functions in the
  space for which $D^\downarrow$ is a reproducing kernel.
\end{itemize}
Note that experience with Slepian functions shows that the eigenvalues
often separate into a set of values close to $1$ and some others which
are almost $0$. This also demonstrates the difficulty of finding a
numerically stable kernel $D^\downarrow$, since then some $\tau_k^{-2}$ are
very large.

\begin{rem}
Since the scaling functions $\Phi_J$ provide us with different
approximations $\Phi_J\ast F$ to $F$, it also appears to be useful to
look at differences $\Psi_J\coloneqq \Phi_{J+1}-\Phi_J$ such that
\begin{equation}
  \Phi_{J+1}\ast G =\Phi_J\ast G +\Psi_J\ast G\,.
\end{equation}
Here, $\Psi_J\ast G$ can be regarded as the detail information with is
added to the approximation $\Phi_J\ast G$ at scale $J$ to obtain the
approximation $\Phi_{J+1}\ast G$ at the next scale.  In analogy to
common wavelet theories, where such scale-step properties also exist,
the kernels $\Psi_J$ can be called wavelets here.
\end{rem}

\section{Some numerical tests}\label{Sec:numerics}
This paper generalizes an approach presented
in \cite{PlattnerSimons2017} for the downward continuation of
geophysically relevant potentials. Their application has served as a
thread in this paper to show that the general setup, indeed, includes
this particular example. Rather than experimenting with the same
examples again, we demonstrate the applicability of the general Ansatz
to other inverse problems by discussing some enlightening problems on
the 1-sphere. All numerical calculations were done
with \texttt{MatlabR2015b}.

\subsection{Identity}

We start with an approximation problem. The Hilbert spaces and the
operator (which is the identity operator for an approximation problem)
are chosen as follows:
\begin{align}
\cX&\coloneqq\rL^2[0,2\pi],&
\cY&\coloneqq\cX,&
\cZ&\coloneqq\rL^2\left[\frac{\pi}{2},\frac{3\pi}{2}\right],\\
T&\coloneqq \mathrm{Id},&
\sigma_n&\coloneqq1 \quad \text{for all }n\,.
\end{align}
Note that $\rL^2[0,2\pi]$ is isometric and isomorphic to
$\rL^2(\mathbb{S}^1)$, where $\mathbb{S}^1$ is the 1-sphere.

Moreover, we need orthonormal basis systems for the Hilbert spaces
involved. We take here a common system, for $x\in[0,2\pi]$,
\begin{align}
u_{0,1}(x)&\coloneqq \frac{1}{\sqrt{2\pi}},\\
u_{n,1}(x)&\coloneqq \frac{1}{\sqrt{\pi}}\, \cos(nx),& 
u_{n,2}(x)&\coloneqq \frac{1}{\sqrt{\pi}}\, \sin(nx),\quad n\geq 1\,,\\
v_{n,j}(x)&\coloneqq u_{n,j}(x)\quad \text{for all }n,j\,.
\end{align}
For calculating the Slepian functions, the bandlimit is set to $N\coloneqq 50$. Moreover, we use the transformation 
\begin{equation}
k(n,j)=\begin{cases}
 2(n-1)+j+1, & \text{if } n\geq 1 \\ 
1, &\text{if } n=0
\end{cases}\,
\end{equation}
to have a single index only.  A selection of the Slepian functions on
the 1-sphere $\mathbb{S}^1$ with largest and lowest eigenvalues is
shown in Figure \ref{Fig:Slepian_Approx}. It can be seen that the set
of Slepian functions can be subdivided into functions with a strong
localization in $R=[0.5\pi,1.5\pi]$ and other functions which
concentrate on the complement $D\setminus
R=[0,0.5\pi[\,\cup\,]1.5\pi,2\pi]$. This is also confirmed by the
eigenvalues, which are shown in Figure \ref{Fig:eigenvalues}. For
numerical reasons, we only consider Slepian functions for which
$\tau_k\geq 0.1\%\cdot\tau_1$ in all our calculations.

\begin{figure}
\centering
\input{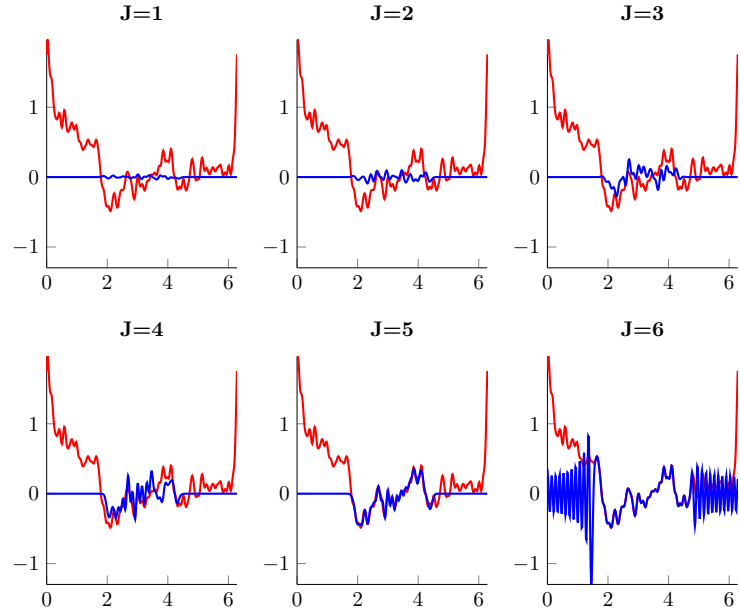}
\caption{Solution (red) and multi-scale approximations (blue) at different scales $J$, shown on the whole domain $D$: at a sufficiently large scale, we obtain a very good approximation to the projection $\cP F$ of $F$ to the subinterval $R$.}
\label{Fig:Approx_sol_on_D}
\end{figure}

\begin{figure}
\centering
\input{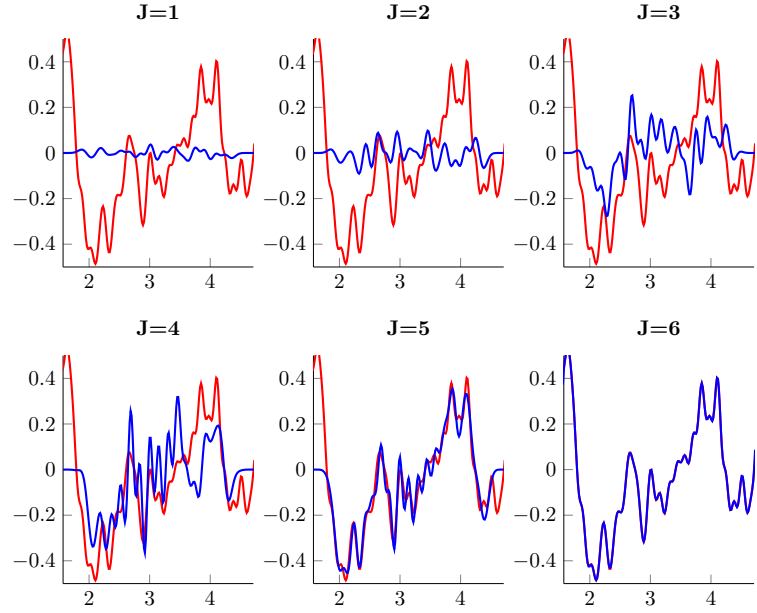}
\caption{Solution (red) and multi-scale approximations (blue) at different scales $J$, shown on the subdomain $R$: at scale $J=6$, the projection $\cP F$ can hardly be distinguished from the approximation $\Phi_J\ast G$.}
\label{Fig:Approx_sol_on_R}
\end{figure}

The Fourier coefficients of the contrived solution $F$ are chosen by
\begin{equation}
\xip{F}{u_k} \coloneqq \left(1+\varepsilon_k\right)\,\frac{1}{k}, \quad k=1,\dots,2N+1\,.
\end{equation}
The $\varepsilon_k$ are standard normally distributed random
numbers. The corresponding function is represented by the red graphs
in Figures \ref{Fig:Approx_sol_on_D} and \ref{Fig:Approx_sol_on_R}.
The right-hand side is calculated as $G=\cP\cT
F=\left.F\right|_{[0.5\pi,1.5\pi]}$ on an equidistant grid of 1001
points in $[0.5\pi,1.5\pi]$. This right-hand side is contaminated with
noise by replacing $G$ with $G+0.01\,\tilde{\varepsilon}_k$ (see
Figure \ref{Fig:Appr_rhs}), where the $\tilde{\varepsilon}_k$ are
standard normally distributed random variables (${\varepsilon}_k$ and
$\tilde{\varepsilon}_k$ were obtained with the \texttt{MATLAB}
function \texttt{randn}). Moreover, the functions $\varphi_J$ are
chosen as the generators of the Shannon scaling function (see
e.g.\ \cite{OxfordBuch}) such that
\begin{equation}\label{Eq:Shannonscf}
  \varphi_J(k) = \begin{cases}
  1, & \text{if } k< 2^J\\
  0, & \text{else}
  \end{cases}
  , \quad J,k\in\N_0\,.
\end{equation}
For the convolution $\zip{\Phi_J(x,\cdot)}{G}$, a composite Simpson's rule was used. The points $x_i$ used for plotting $\zip{\Phi_J(x,\cdot)}{G}$ are on an equidistant grid of 401 points in $[0,2\pi]$.

The root mean square error $(\frac{1}{M}\,\sum_{i=1}^M (F(x_i)-(\Phi_J\ast G)(x_i))^2)^{1/2}$ is calculated only for points $x_i$ in $R=[0.5\pi,1.5\pi]$ and is shown in Table \ref{Tab:Id}. The approximation error clearly decreases and then stagnates at a low level (note that the truncation condition $\tau_k<0.1\%\cdot\tau_1$ is achieved in this example for $k=60$; hence, we have here that $\varphi_{J_1}(k)=\varphi_{J_2}(k)$ for all $k=1,\dots,101$, if $J_1,J_2\geq 6$). Note that the values of $F$ vary within $R$ between $-0.5$ and $0.5$. The obtained approximations are shown in Figures \ref{Fig:Approx_sol_on_D} and \ref{Fig:Approx_sol_on_R}. We can see that the chosen function $F$ is well approximated on the interval $R$. For the larger scale $J=6$, some boundary effects\footnote{We experienced in our experiments that a finer quadrature grid of $10,001$ points for the Simpson rule reduces these effects in their amplitude such that they can partially also occur due to inaccuracies in the numerical integration; however, also with this finer grid, the effects were still clearly visible.} occur, which shows that, in some cases, smoother approximations at lower scales (like here for $J=5$), which are still close to the exact solution but do not show such boundary effects, might be preferred.
\begin{figure}
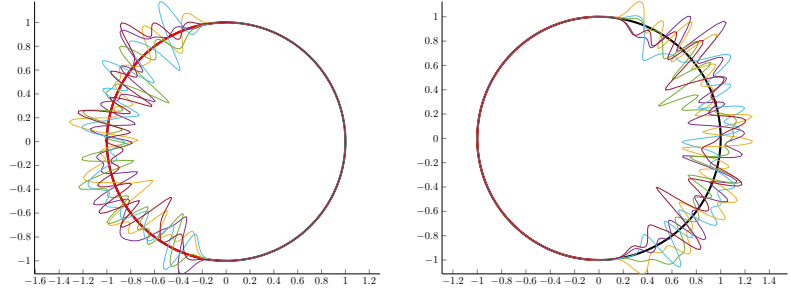

\centering
\input{Fig_Appr_Slep_large_ev.tex}
\hspace*{0.2cm}
\input{Fig_Appr_Slep_small_ev.tex}
\caption{Eigenfunctions $g_1,\dots,g_6$ corresponding to the largest eigenvalues and $g_{96},\dots,g_{101}$ corresponding to the smallest eigenvalues, the subdomain $R$ is shown in red; note that each Slepian function was multiplied with the same factor to scale the amplitudes for better visibility.}
\label{Fig:Slepian_Approx}
\end{figure}
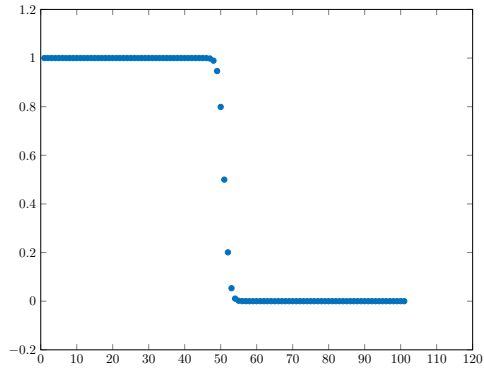
\begin{figure}
\centering
%
\definecolor{mycolor1}{rgb}{0.00000,0.44700,0.74100}%
\begin{tikzpicture}[scale=0.5]

\begin{axis}[%
width=4.521in,
height=3.566in,
at={(0.758in,0.481in)},
scale only axis,
xmin=0,
xmax=120,
ymin=-0.2,
ymax=1.2,
axis background/.style={fill=white}
]
\addplot [color=mycolor1,only marks,mark=*,mark options={solid},forget plot]
  table[row sep=crcr]{%
1	1\\
2	1\\
3	1\\
4	1\\
5	1\\
6	1\\
7	1\\
8	1\\
9	1\\
10	1\\
11	1\\
12	1\\
13	1\\
14	1\\
15	1\\
16	1\\
17	1\\
18	1\\
19	1\\
20	1\\
21	1\\
22	1\\
23	1\\
24	1\\
25	1\\
26	1\\
27	0.999999999999999\\
28	0.999999999999999\\
29	0.999999999999999\\
30	0.999999999999999\\
31	0.999999999999999\\
32	0.999999999999999\\
33	0.999999999999999\\
34	0.999999999999999\\
35	0.999999999999998\\
36	0.99999999999998\\
37	0.999999999999727\\
38	0.99999999999638\\
39	0.999999999955136\\
40	0.999999999479589\\
41	0.999999994360791\\
42	0.999999943039844\\
43	0.999999465180979\\
44	0.999995348653202\\
45	0.999962707009083\\
46	0.999726202359545\\
47	0.998179169837806\\
48	0.989265172611167\\
49	0.946644484830518\\
50	0.799061458634963\\
51	0.499999999999999\\
52	0.200938541365035\\
53	0.0533555151694834\\
54	0.0107348273888352\\
55	0.00182083016219454\\
56	0.000273797640452691\\
57	3.72929909153044e-05\\
58	4.65134679846005e-06\\
59	5.34819020715571e-07\\
60	5.69601550839168e-08\\
61	5.63921044894273e-09\\
62	5.20411637216494e-10\\
63	4.48646760262125e-11\\
64	3.62009167003476e-12\\
65	2.73598629715599e-13\\
66	1.91827986890296e-14\\
67	1.47692185151901e-15\\
68	-1.07408402170279e-15\\
69	-1.02489292856735e-15\\
70	9.16040896473652e-16\\
71	-8.95445334715226e-16\\
72	-8.60522620888916e-16\\
73	7.61655847773018e-16\\
74	-6.78572032578671e-16\\
75	-6.75526553296299e-16\\
76	6.68383438643344e-16\\
77	-6.64180810058828e-16\\
78	6.13688558007415e-16\\
79	5.72434859243452e-16\\
80	-4.80475524838579e-16\\
81	4.42107538362402e-16\\
82	4.27887520844726e-16\\
83	3.88876160084776e-16\\
84	-3.86528301324784e-16\\
85	-3.62180093994903e-16\\
86	-3.23071883302424e-16\\
87	-3.18853748522749e-16\\
88	-3.01140963894808e-16\\
89	2.91097805445652e-16\\
90	2.62860131698656e-16\\
91	2.59681608621543e-16\\
92	-2.30204954802363e-16\\
93	-1.9685383333837e-16\\
94	-1.82278362084402e-16\\
95	-8.32579211035325e-17\\
96	8.1903587514284e-17\\
97	-6.93816009196104e-17\\
98	5.30558679417355e-17\\
99	-2.90790467821559e-17\\
100	-1.07056974886127e-17\\
101	-1.00125218763979e-17\\
};
\end{axis}
\end{tikzpicture}%
\caption{Eigenvalues (sorted) for the Slepian localization problem: there is a sharp transition from strong to weak localization.}
\label{Fig:eigenvalues}
\end{figure}
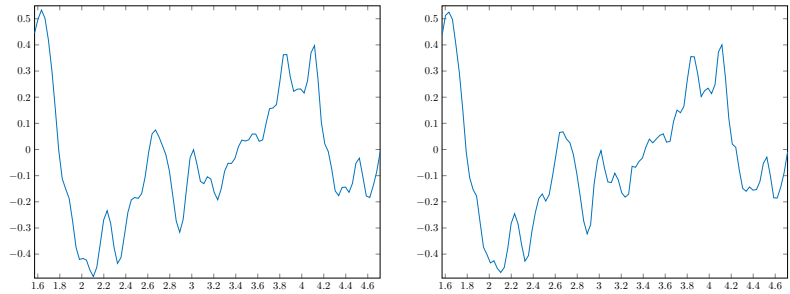
\begin{figure}
\centering
%
\definecolor{mycolor1}{rgb}{0.00000,0.44700,0.74100}%
\begin{tikzpicture}[scale=0.4]

\begin{axis}[%
width=4.521in,
height=3.566in,
at={(0.758in,0.481in)},
scale only axis,
xmin=1.5707963267949,
xmax=4.71238898038469,
ymin=-0.491696690922447,
ymax=0.549487668109766,
axis background/.style={fill=white}
]
\addplot [color=mycolor1,solid,forget plot]
  table[row sep=crcr]{%
1.5707963267949	0.445557131906617\\
1.60221225333079	0.49912038987572\\
1.63362817986669	0.533157637228094\\
1.66504410640259	0.502071698831883\\
1.69646003293849	0.418543816521522\\
1.72787595947439	0.303128477983481\\
1.75929188601028	0.154154687608098\\
1.79070781254618	-0.00425344857846038\\
1.82212373908208	-0.111441185583064\\
1.85353966561798	-0.150890693655581\\
1.88495559215388	-0.186893595413932\\
1.91637151868977	-0.273637800543758\\
1.94778744522567	-0.374134815352615\\
1.97920337176157	-0.420478425760291\\
2.01061929829747	-0.416441215821185\\
2.04203522483337	-0.422662754474024\\
2.07345115136926	-0.460745846165529\\
2.10486707790516	-0.485767410731744\\
2.13628300444106	-0.450587669624224\\
2.16769893097696	-0.361058367492127\\
2.19911485751286	-0.269100851415512\\
2.23053078404875	-0.234063491258637\\
2.26194671058465	-0.281210269376273\\
2.29336263712055	-0.374193631790715\\
2.32477856365645	-0.435348673837101\\
2.35619449019234	-0.413447629144321\\
2.38761041672824	-0.328350223791194\\
2.41902634326414	-0.240861526412378\\
2.45044226980004	-0.192465268171267\\
2.48185819633594	-0.183288328511774\\
2.51327412287183	-0.186902201710074\\
2.54469004940773	-0.168104387670453\\
2.57610597594363	-0.104274651385591\\
2.60752190247953	-0.0117131584113877\\
2.63893782901543	0.0595537575577623\\
2.67035375555132	0.0747433000552827\\
2.70176968208722	0.0483357052844319\\
2.73318560862312	0.0149273594642415\\
2.76460153515902	-0.0202181494208772\\
2.79601746169492	-0.0809732425637255\\
2.82743338823081	-0.176330899118367\\
2.85884931476671	-0.27352567336351\\
2.89026524130261	-0.316234356371631\\
2.92168116783851	-0.26744089599113\\
2.95309709437441	-0.14685672522771\\
2.9845130209103	-0.031879184113178\\
3.0159289474462	-0.00081105507981136\\
3.0473448739821	-0.0566513004330042\\
3.078760800518	-0.122095043465106\\
3.1101767270539	-0.130414680946866\\
3.14159265358979	-0.104177152161199\\
3.17300858012569	-0.112661280666576\\
3.20442450666159	-0.16463218757659\\
3.23584043319749	-0.19223272565392\\
3.26725635973339	-0.150628723079789\\
3.29867228626928	-0.082864312626769\\
3.33008821280518	-0.0526272955605252\\
3.36150413934108	-0.0531681063090664\\
3.39292006587698	-0.0335047113454148\\
3.42433599241287	0.010773903379156\\
3.45575191894877	0.0360544217767915\\
3.48716784548467	0.0323164993621395\\
3.51858377202057	0.037263882684545\\
3.54999969855647	0.0593841634186861\\
3.58141562509236	0.0588576970989297\\
3.61283155162826	0.0313291660349338\\
3.64424747816416	0.0366575861864791\\
3.67566340470006	0.0997235979674015\\
3.70707933123596	0.155936502878423\\
3.73849525777185	0.159059783107529\\
3.76991118430775	0.171811617802926\\
3.80132711084365	0.260926171434164\\
3.83274303737955	0.362343850276601\\
3.86415896391545	0.362819633195279\\
3.89557489045134	0.277516279066712\\
3.92699081698724	0.22248752516219\\
3.95840674352314	0.23060420759078\\
3.98982267005904	0.231674875453968\\
4.02123859659494	0.215950629556445\\
4.05265452313083	0.263200350586892\\
4.08407044966673	0.369738001068723\\
4.11548637620263	0.397396719007314\\
4.14690230273853	0.271716232488059\\
4.17831822927442	0.10357375470854\\
4.20973415581032	0.0209166284307319\\
4.24115008234622	-0.00782725043740071\\
4.27256600888212	-0.0748394656454567\\
4.30398193541802	-0.15819782840069\\
4.33539786195391	-0.176267713056059\\
4.36681378848981	-0.145202742519392\\
4.39822971502571	-0.143863286948524\\
4.42964564156161	-0.16367735383749\\
4.46106156809751	-0.129463398534666\\
4.4924774946334	-0.0524655914970962\\
4.5238934211693	-0.0330149981024505\\
4.5553093477052	-0.103101152317435\\
4.5867252742411	-0.177738600963011\\
4.618141200777	-0.183356643841049\\
4.64955712731289	-0.139584629061801\\
4.68097305384879	-0.0858956816673784\\
4.71238898038469	-0.0108864005504366\\
};
\end{axis}
\end{tikzpicture}%
\hspace*{0.2cm}
%
\definecolor{mycolor1}{rgb}{0.00000,0.44700,0.74100}%
\begin{tikzpicture}[scale=0.4]

\begin{axis}[%
width=4.521in,
height=3.566in,
at={(0.758in,0.481in)},
scale only axis,
xmin=1.5707963267949,
xmax=4.71238898038469,
ymin=-0.491696690922447,
ymax=0.549487668109766,
axis background/.style={fill=white}
]
\addplot [color=mycolor1,solid,forget plot]
  table[row sep=crcr]{%
1.5707963267949	0.435851888310912\\
1.60221225333079	0.512817508352385\\
1.63362817986669	0.525352907625298\\
1.66504410640259	0.497350022339904\\
1.69646003293849	0.401279746063234\\
1.72787595947439	0.294389125862974\\
1.75929188601028	0.151917358314631\\
1.79070781254618	-0.0102388239870804\\
1.82212373908208	-0.107941651458974\\
1.85353966561798	-0.153592937053444\\
1.88495559215388	-0.178198681970014\\
1.91637151868977	-0.277403720591365\\
1.94778744522567	-0.373936718174339\\
1.97920337176157	-0.401460876556267\\
2.01061929829747	-0.433722680967019\\
2.04203522483337	-0.425225812693952\\
2.07345115136926	-0.454645507480453\\
2.10486707790516	-0.470193556492258\\
2.13628300444106	-0.450931259855888\\
2.16769893097696	-0.379570426571253\\
2.19911485751286	-0.282792803763509\\
2.23053078404875	-0.245510908069062\\
2.26194671058465	-0.284909805608592\\
2.29336263712055	-0.362437687217645\\
2.32477856365645	-0.426519483385336\\
2.35619449019234	-0.405833246627557\\
2.38761041672824	-0.313806502583721\\
2.41902634326414	-0.241114995066535\\
2.45044226980004	-0.1871326726587\\
2.48185819633594	-0.170078429695208\\
2.51327412287183	-0.197280673767529\\
2.54469004940773	-0.173134856467347\\
2.57610597594363	-0.0995221187150254\\
2.60752190247953	-0.0178979195956138\\
2.63893782901543	0.064862777053896\\
2.67035375555132	0.0678006824804781\\
2.70176968208722	0.0406491676547607\\
2.73318560862312	0.0256915259353599\\
2.76460153515902	-0.0187416292630927\\
2.79601746169492	-0.0923085884572339\\
2.82743338823081	-0.180339176099223\\
2.85884931476671	-0.273574398260114\\
2.89026524130261	-0.322221985421819\\
2.92168116783851	-0.287221389319651\\
2.95309709437441	-0.133422580181546\\
2.9845130209103	-0.0412909728500985\\
3.0159289474462	-0.00416472801842999\\
3.0473448739821	-0.0707299964423455\\
3.078760800518	-0.123950039077105\\
3.1101767270539	-0.12622444819109\\
3.14159265358979	-0.090643141877748\\
3.17300858012569	-0.116264310537492\\
3.20442450666159	-0.165387946112705\\
3.23584043319749	-0.182057686511275\\
3.26725635973339	-0.170605072483459\\
3.29867228626928	-0.0642494060188033\\
3.33008821280518	-0.0680394651400735\\
3.36150413934108	-0.0458253455814288\\
3.39292006587698	-0.0325873699892269\\
3.42433599241287	0.00955130374014521\\
3.45575191894877	0.0396293460854333\\
3.48716784548467	0.0257023419065568\\
3.51858377202057	0.0402519304844408\\
3.54999969855647	0.0538305063963294\\
3.58141562509236	0.0601301474644378\\
3.61283155162826	0.0279356248924158\\
3.64424747816416	0.0317960867957991\\
3.67566340470006	0.107142866900832\\
3.70707933123596	0.15046719657264\\
3.73849525777185	0.140712173427704\\
3.76991118430775	0.165199546696155\\
3.80132711084365	0.268088219488655\\
3.83274303737955	0.355159271666678\\
3.86415896391545	0.354283684506309\\
3.89557489045134	0.288952031488698\\
3.92699081698724	0.202720337819697\\
3.95840674352314	0.224600007574811\\
3.98982267005904	0.234376142795958\\
4.02123859659494	0.213793044051385\\
4.05265452313083	0.247884671763769\\
4.08407044966673	0.373203249183495\\
4.11548637620263	0.40014233483686\\
4.14690230273853	0.277139102440439\\
4.17831822927442	0.11830794456911\\
4.20973415581032	0.0208112383587769\\
4.24115008234622	0.00862305762657309\\
4.27256600888212	-0.0761945024744159\\
4.30398193541802	-0.148236776198539\\
4.33539786195391	-0.159699363275411\\
4.36681378848981	-0.143263479793928\\
4.39822971502571	-0.155423617012997\\
4.42964564156161	-0.153232101244718\\
4.46106156809751	-0.122393850844978\\
4.4924774946334	-0.0525336805191489\\
4.5238934211693	-0.0293937133262555\\
4.5553093477052	-0.0980755629705495\\
4.5867252742411	-0.184812758178826\\
4.618141200777	-0.185345370366977\\
4.64955712731289	-0.144832434088117\\
4.68097305384879	-0.088998561691535\\
4.71238898038469	-0.00910780830368501\\
};
\end{axis}
\end{tikzpicture}%
\caption{Right-hand side $G$ for the approximation problem without noise (left) and after adding the noise (right)}
\label{Fig:Appr_rhs}
\end{figure}
\begin{table}
\centering
\begin{tabular}{|l|l|}
scale & error\\
1 & 0.23846\\
2 & 0.23566\\
3 & 0.21851\\
4 & 0.18388\\
5 & 0.12127\\
6 & 0.0024189\\
7 & 0.0024189\\
\end{tabular}
\caption{Errors (rms) for the pure approximation problem depending on the scale $J$ for the Shannon scaling function $\Phi_J$.}
\label{Tab:Id}
\end{table}
\clearpage
\subsection{An inverse problem}
We now consider an ill-posed inverse problem. The spaces $\cX$, $\cY$, and $\cZ$ as well as their orthonormal basis systems are chosen like above. However, the singular values are now given by
\begin{equation}
\sigma_{n,j}=\frac{1}{n+1} \quad \text{for all }n,j\,.\\
\end{equation}
In Figures \ref{Fig:IP_Slepian_ip} and \ref{Fig:IP_eigenvalues}, we can see that the Slepian functions and their eigenvalues are indeed influenced by the ill-posed nature of the problem.

Again, the bandlimit is set to $N\coloneqq 50 $. We also take the same function $F$ as the solution of $\cP\cT F=G$. The right-hand side $G$ is shown in Figure \ref{Fig:IP_rhs}.

\begin{figure}
\centering
\input{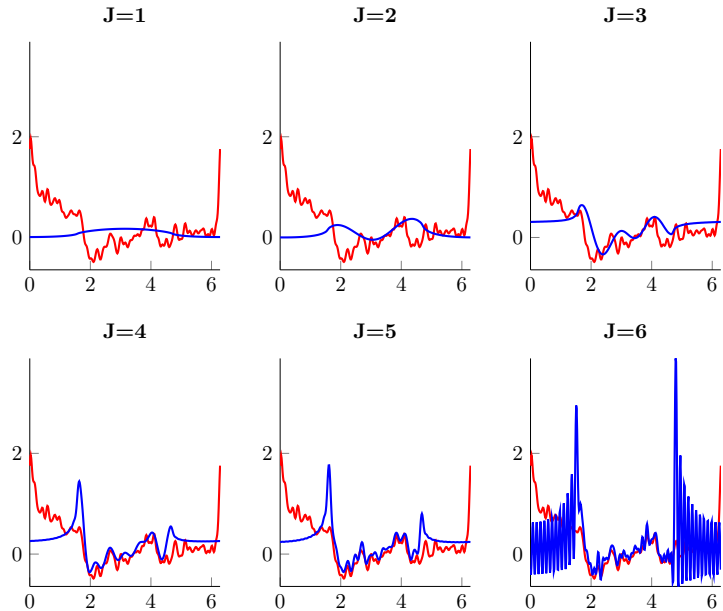}
\caption{Solution of the inverse problem (red) and multi-scale approximations (blue) at different scales $J$, shown on the whole domain $D$: in view of the ill-posed nature of the problem, the approximations are still rather close to the solution $F$ on the subdomain $R$.}
\label{Fig:ip_sol_on_D}
\end{figure}

\begin{figure}
\centering
\input{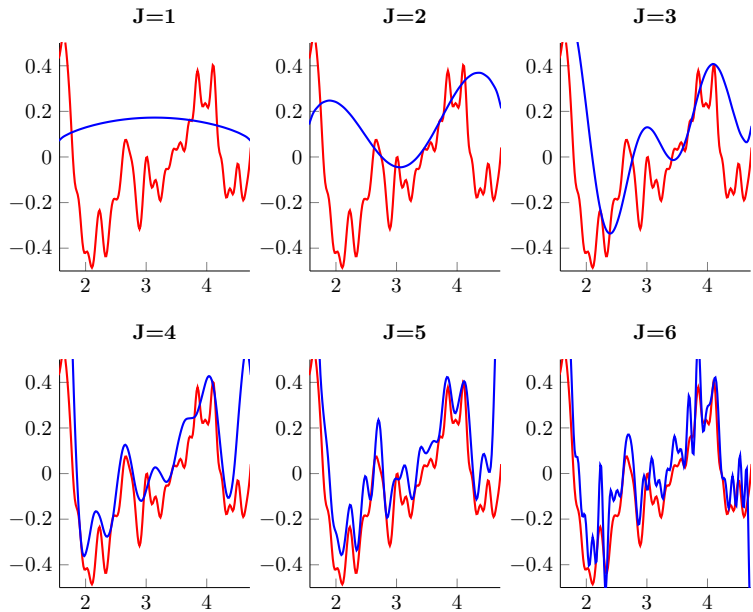}
\caption{Solution of the inverse problem (red) and multi-scale approximations (blue) at different scales $J$, shown on the subdomain $R$: the approximations are relatively close to the exact solution. Depending on the scale, the approximations are more or less filtered as the varying smoothness shows.}
\label{Fig:ip_sol_on_R}
\end{figure}

The rest of the numerical calculations is performed like above.
The results are shown in Table \ref{Tab:invprobl} and Figures \ref{Fig:ip_sol_on_D} and \ref{Fig:ip_sol_on_R} (here, the truncation condition $\tau_k<0.1\%\cdot\tau_1$ is reached for $k=58$ such that the approximations again stagnate from scale $J=6$). Clearly, the noise has much more influence on the solution of the ill-posed problem. However, the approximations at sufficiently large scales are still rather close to the exact (noise-free) solution.

We can also see that the multiscale approach is appropriate for smoothing the solution. For example, scales $J=3$ and $J=4$ reveal trends in the solution which are smooth and coarse (i.e.\ associated to a low frequency). This is, for example, useful, if a very noisy signal can be expected or if one is interested in separating the phenomena of different \lq wavelengths\rq\ (in a more abstract sense) in the solution.
\begin{figure}
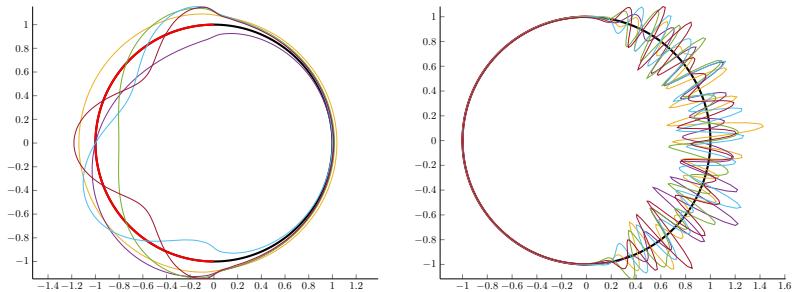

\centering
\input{Fig_ip_Slep_large_ev.tex}
\hspace*{0.2cm}
\input{Fig_ip_Slep_small_ev.tex}
\caption{Eigenfunctions $g_1,\dots,g_6$ and $g_{96},\dots,g_{101}$ for the case of the chosen inverse problem, the subdomain $R$ is shown in red; note that each Slepian function was multiplied with the same factor to scale the amplitudes for better visibility.}
\label{Fig:IP_Slepian_ip}
\end{figure}
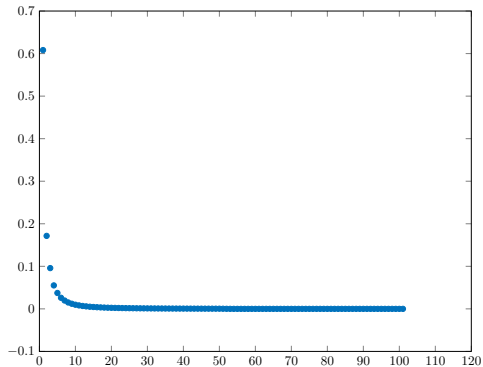
\begin{figure}
\centering
%
\definecolor{mycolor1}{rgb}{0.00000,0.44700,0.74100}%
\begin{tikzpicture}[scale=0.5]

\begin{axis}[%
width=4.521in,
height=3.566in,
at={(0.758in,0.481in)},
scale only axis,
xmin=0,
xmax=120,
ymin=-0.1,
ymax=0.7,
axis background/.style={fill=white}
]
\addplot [color=mycolor1,only marks,mark=*,mark options={solid},forget plot]
  table[row sep=crcr]{%
1	0.608226469661418\\
2	0.171535795075713\\
3	0.0958151916894951\\
4	0.0551252607008364\\
5	0.0373550875272519\\
6	0.0260661941339458\\
7	0.0196088827013379\\
8	0.0150286091192586\\
9	0.0120230311455971\\
10	0.00973831132332224\\
11	0.00810885183097715\\
12	0.00681116230101405\\
13	0.00583243237322443\\
14	0.00502664272639181\\
15	0.00439402349800464\\
16	0.00386010119169139\\
17	0.00342806091009527\\
18	0.00305636188303481\\
19	0.00274840559064795\\
20	0.00247939335916996\\
21	0.00225224411372462\\
22	0.00205135430754372\\
23	0.00187905517814722\\
24	0.00172511624745875\\
25	0.00159134072224757\\
26	0.00147080200752512\\
27	0.00136486786367973\\
28	0.00126873217629842\\
29	0.00118341378659907\\
30	0.0011055183287892\\
31	0.00103578713186379\\
32	0.000971795862068726\\
33	0.000914065056238882\\
34	0.000860854842515931\\
35	0.000812510264303032\\
36	0.000767784004831091\\
37	0.000726883714702431\\
38	0.000688918674771636\\
39	0.000653995367458835\\
40	0.000621473785709315\\
41	0.000591401808389432\\
42	0.000563290716875705\\
43	0.00053719501613602\\
44	0.00051264916024872\\
45	0.000489843015782084\\
46	0.000468089267022527\\
47	0.000448036142357169\\
48	0.000428046951671927\\
49	0.000410366692933139\\
50	0.000383823476039615\\
51	0.000285219335057138\\
52	0.000133451130076346\\
53	4.0296987822666e-05\\
54	8.81092283785266e-06\\
55	1.604679659739e-06\\
56	2.52165449309656e-07\\
57	3.61936549683202e-08\\
58	4.66281640287455e-09\\
59	5.61787342477232e-10\\
60	6.14777267637112e-11\\
61	6.36315733426638e-12\\
62	6.0148964242503e-13\\
63	5.41702545477188e-14\\
64	4.46731904250959e-15\\
65	3.53258159007844e-16\\
66	-1.98673953573473e-17\\
67	1.9301740735991e-17\\
68	8.84415730645235e-18\\
69	5.56576738684201e-18\\
70	-5.54735298520839e-18\\
71	-4.66344988458859e-18\\
72	4.25608744799052e-18\\
73	3.79400550688197e-18\\
74	-3.566879118707e-18\\
75	3.37308523956127e-18\\
76	-3.04195812552155e-18\\
77	-2.2512098345413e-18\\
78	-2.06987150091599e-18\\
79	2.04554709083252e-18\\
80	-2.02119930537475e-18\\
81	1.66312093575082e-18\\
82	-1.51879338716523e-18\\
83	1.25555255960584e-18\\
84	-1.17757838544211e-18\\
85	1.04911070317375e-18\\
86	1.0029607841699e-18\\
87	-9.4645062373878e-19\\
88	7.77476749868696e-19\\
89	-6.19644908291185e-19\\
90	6.06076856906324e-19\\
91	5.31573637110239e-19\\
92	-5.17748252016505e-19\\
93	4.99012390642083e-19\\
94	-4.71126910108912e-19\\
95	-4.45981304182265e-19\\
96	-4.07965369625111e-19\\
97	-3.04658399916936e-19\\
98	-2.3781540829297e-19\\
99	1.79806977182637e-19\\
100	1.46084739121366e-19\\
101	7.16157115802192e-20\\
};
\end{axis}
\end{tikzpicture}%
\caption{Eigenvalues (sorted) for the Slepian localization problem and the case of the chosen inverse problem: the ill-posedness is also reflected in the eigenvalues.}
\label{Fig:IP_eigenvalues}
\end{figure}
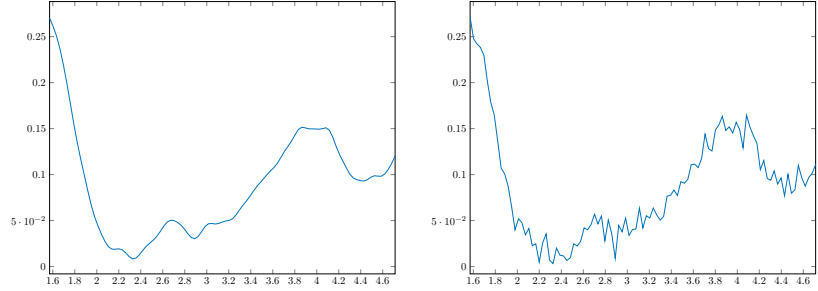
\begin{figure}
\centering
%
\definecolor{mycolor1}{rgb}{0.00000,0.44700,0.74100}%
\begin{tikzpicture}[scale=0.4]

\begin{axis}[%
width=4.521in,
height=3.566in,
at={(0.758in,0.481in)},
scale only axis,
xmin=1.5707963267949,
xmax=4.71238898038469,
ymin=-0.00793970297215048,
ymax=0.288134496563749,
axis background/.style={fill=white}
]
\addplot [color=mycolor1,solid,forget plot]
  table[row sep=crcr]{%
1.5707963267949	0.270126904936695\\
1.60221225333079	0.261307149990951\\
1.63362817986669	0.250666602081418\\
1.66504410640259	0.23666650606244\\
1.69646003293849	0.219294031728668\\
1.72787595947439	0.199275167427915\\
1.75929188601028	0.177425958227946\\
1.79070781254618	0.1552868823332\\
1.82212373908208	0.134904316754876\\
1.85353966561798	0.116959308142935\\
1.88495559215388	0.100145658032596\\
1.91637151868977	0.0833577065553755\\
1.94778744522567	0.0675618712225272\\
1.97920337176157	0.0544660555586852\\
2.01061929829747	0.0441602984434179\\
2.04203522483337	0.03534783027788\\
2.07345115136926	0.0275261809831911\\
2.10486707790516	0.021725510994913\\
2.13628300444106	0.0189863168009639\\
2.16769893097696	0.0188211316632043\\
2.19911485751286	0.0192911395726288\\
2.23053078404875	0.0182954032325753\\
2.26194671058465	0.0150633171988052\\
2.29336263712055	0.0109072244648672\\
2.32477856365645	0.0083814658879128\\
2.35619449019234	0.00921770253127092\\
2.38761041672824	0.0129866666692519\\
2.41902634326414	0.0178117000374619\\
2.45044226980004	0.0221027568057798\\
2.48185819633594	0.0255056214815992\\
2.51327412287183	0.0287028013424705\\
2.54469004940773	0.032710770539994\\
2.57610597594363	0.0379949568648165\\
2.60752190247953	0.0437657472754206\\
2.63893782901543	0.0483045419032789\\
2.67035375555132	0.0503449220845223\\
2.70176968208722	0.0499995445939667\\
2.73318560862312	0.0481216358915116\\
2.76460153515902	0.0451191987311378\\
2.79601746169492	0.040899937944642\\
2.82743338823081	0.0358921710043355\\
2.85884931476671	0.0316555956177714\\
2.89026524130261	0.0302248497259758\\
2.92168116783851	0.0326796167046127\\
2.95309709437441	0.0380408360628624\\
2.9845130209103	0.0435446995600437\\
3.0159289474462	0.04657496956357\\
3.0473448739821	0.0468188465273892\\
3.078760800518	0.0462929452390888\\
3.1101767270539	0.0469143942902741\\
3.14159265358979	0.0484701312764642\\
3.17300858012569	0.0495834406509773\\
3.20442450666159	0.0502279746294398\\
3.23584043319749	0.0520996238103516\\
3.26725635973339	0.0562942386519212\\
3.29867228626928	0.0617780931363191\\
3.33008821280518	0.0668993796595308\\
3.36150413934108	0.0715206476016594\\
3.39292006587698	0.0765744571707093\\
3.42433599241287	0.0821124007010499\\
3.45575191894877	0.0871863999864438\\
3.48716784548467	0.0915138828559978\\
3.51858377202057	0.0958070772408075\\
3.54999969855647	0.100266926679862\\
3.58141562509236	0.104252739356639\\
3.61283155162826	0.107937057276555\\
3.64424747816416	0.112705516384746\\
3.67566340470006	0.11901493377902\\
3.70707933123596	0.125388702893536\\
3.73849525777185	0.130772511202773\\
3.76991118430775	0.136313689889696\\
3.80132711084365	0.143011989447506\\
3.83274303737955	0.148985869979236\\
3.86415896391545	0.151477479671504\\
3.89557489045134	0.150802886080007\\
3.92699081698724	0.149783497844097\\
3.95840674352314	0.149706465535989\\
3.98982267005904	0.149599596384445\\
4.02123859659494	0.149325993261449\\
4.05265452313083	0.150055874590624\\
4.08407044966673	0.150870650792486\\
4.11548637620263	0.148181197414046\\
4.14690230273853	0.140302065444581\\
4.17831822927442	0.130094607839123\\
4.20973415581032	0.121176565389632\\
4.24115008234622	0.113780575500832\\
4.27256600888212	0.106508749907977\\
4.30398193541802	0.100004747097394\\
4.33539786195391	0.0960460203517868\\
4.36681378848981	0.0943801918926293\\
4.39822971502571	0.0933711315973357\\
4.42964564156161	0.0930159241395353\\
4.46106156809751	0.0945248151226381\\
4.4924774946334	0.0972404520771842\\
4.5238934211693	0.0987409916304963\\
4.5553093477052	0.0983365860523628\\
4.5867252742411	0.0982187273138022\\
4.618141200777	0.100598957283563\\
4.64955712731289	0.105546210551051\\
4.68097305384879	0.112188677543171\\
4.71238898038469	0.120229220954367\\
};
\end{axis}
\end{tikzpicture}%
\hspace*{0.2cm}
%
\definecolor{mycolor1}{rgb}{0.00000,0.44700,0.74100}%
\begin{tikzpicture}[scale=0.4]

\begin{axis}[%
width=4.521in,
height=3.566in,
at={(0.758in,0.481in)},
scale only axis,
xmin=1.5707963267949,
xmax=4.71238898038469,
ymin=-0.00793970297215048,
ymax=0.288134496563749,
axis background/.style={fill=white}
]
\addplot [color=mycolor1,solid,forget plot]
  table[row sep=crcr]{%
1.5707963267949	0.271063329130351\\
1.60221225333079	0.247668827035417\\
1.63362817986669	0.242068175847823\\
1.66504410640259	0.238193534535817\\
1.69646003293849	0.229424892774748\\
1.72787595947439	0.20125248579201\\
1.75929188601028	0.178124280313444\\
1.79070781254618	0.164644797567628\\
1.82212373908208	0.135968664513567\\
1.85353966561798	0.106791729080447\\
1.88495559215388	0.10082156835032\\
1.91637151868977	0.0874191552893827\\
1.94778744522567	0.0652823530590766\\
1.97920337176157	0.0400099972295277\\
2.01061929829747	0.0519603197285834\\
2.04203522483337	0.0476561246391113\\
2.07345115136926	0.0344597953114304\\
2.10486707790516	0.0413615523674198\\
2.13628300444106	0.0226022479415869\\
2.16769893097696	0.0247916072530939\\
2.19911485751286	0.004818274597758\\
2.23053078404875	0.0258674473629981\\
2.26194671058465	0.0353532716427788\\
2.29336263712055	0.00691575332527906\\
2.32477856365645	0.00331293291393841\\
2.35619449019234	0.0200778741079271\\
2.38761041672824	0.012416005217772\\
2.41902634326414	0.0116058617377184\\
2.45044226980004	0.00666135324695881\\
2.48185819633594	0.00956115790068498\\
2.51327412287183	0.0247483324472021\\
2.54469004940773	0.0223307915277979\\
2.57610597594363	0.0275089700531969\\
2.60752190247953	0.0420894393696427\\
2.63893782901543	0.0399599803907033\\
2.67035375555132	0.0457474520954974\\
2.70176968208722	0.0568237161374906\\
2.73318560862312	0.0464788778155636\\
2.76460153515902	0.0549758434076469\\
2.79601746169492	0.0280258677886306\\
2.82743338823081	0.0506809187591081\\
2.85884931476671	0.0353316751562149\\
2.89026524130261	0.0089791735671631\\
2.92168116783851	0.0447315667005703\\
2.95309709437441	0.0378744052281223\\
2.9845130209103	0.0528487103790968\\
3.0159289474462	0.0339156658363249\\
3.0473448739821	0.0402543154796701\\
3.078760800518	0.040892168483911\\
3.1101767270539	0.0635490189012036\\
3.14159265358979	0.0421158137265397\\
3.17300858012569	0.0553040618455366\\
3.20442450666159	0.0527814751692966\\
3.23584043319749	0.0636468920264867\\
3.26725635973339	0.0562369160287172\\
3.29867228626928	0.0505887749765147\\
3.33008821280518	0.0543676668148654\\
3.36150413934108	0.0765939444775586\\
3.39292006587698	0.077468918486183\\
3.42433599241287	0.0831614938745254\\
3.45575191894877	0.0772517174763545\\
3.48716784548467	0.0923543551258154\\
3.51858377202057	0.0905882664002728\\
3.54999969855647	0.0947149744918292\\
3.58141562509236	0.110829757074374\\
3.61283155162826	0.111330988605878\\
3.64424747816416	0.107620314220246\\
3.67566340470006	0.117562782671751\\
3.70707933123596	0.143858190071022\\
3.73849525777185	0.128268181469797\\
3.76991118430775	0.125643918565645\\
3.80132711084365	0.148557584484387\\
3.83274303737955	0.15406545611695\\
3.86415896391545	0.163307140402869\\
3.89557489045134	0.14806915980668\\
3.92699081698724	0.151934481235859\\
3.95840674352314	0.14526980600425\\
3.98982267005904	0.15677327306281\\
4.02123859659494	0.149143333357969\\
4.05265452313083	0.129328367321221\\
4.08407044966673	0.164592450290295\\
4.11548637620263	0.151602984092288\\
4.14690230273853	0.142152343607509\\
4.17831822927442	0.13412505116644\\
4.20973415581032	0.10542247412608\\
4.24115008234622	0.115307195625551\\
4.27256600888212	0.095868183303409\\
4.30398193541802	0.0938069902029259\\
4.33539786195391	0.103993955784681\\
4.36681378848981	0.0897327139333266\\
4.39822971502571	0.0962821150426888\\
4.42964564156161	0.077034726701661\\
4.46106156809751	0.100150957758025\\
4.4924774946334	0.0797665748718441\\
4.5238934211693	0.083439537625571\\
4.5553093477052	0.109980063950741\\
4.5867252742411	0.0964159064355459\\
4.618141200777	0.0874352993390279\\
4.64955712731289	0.096791488078162\\
4.68097305384879	0.101432590932489\\
4.71238898038469	0.111035330055991\\
};
\end{axis}
\end{tikzpicture}%
\caption{Right-hand side $G$ for the inverse problem without noise (left) and after adding the noise (right): due to the decreasing singular values, the amplitude of $G$ is smaller than the amplitude of $F$ such that the noise has a stronger influence than in the approximation problem above.}
\label{Fig:IP_rhs}
\end{figure}
\begin{table}
\centering
\begin{tabular}{|l|l|}
scale & error\\
1 & 0.30531\\
2 & 0.32736\\
3 & 0.28118\\
4 & 0.27993\\
5 & 0.28808\\
6 & 0.16788\\
7 & 0.16788
\end{tabular}
\caption{Errors (rms) depending on the scale for the Shannon scaling function for the inverse problem}
\label{Tab:invprobl}
\end{table}

\section{Conclusions}\label{Sec:concl}
We presented a method for the regularization of linear ill-posed
problems as they arise in the geosciences and numerous other
disciplines, where the data are only regionally given, and where the
singular-value decomposition (svd) of the corresponding compact
operator $\cT$ needs to be known only for the global case.  To treat
the case of regional data, we introduced a projection operator $\cP$,
which could be the restriction of functions on the global
domain~$\tilde{D}$ to a regional subdomain~$R$. The idea of the
methodology is based on the interpretation of the quotient of the norm
of the range of $\cP\cT$ and the norm of the preimage as analogous to
the energy ratio as used for the construction of Slepian
functions. The supremum of this quotient is also the operator norm of
$\cP\cT$. Orthonormal ``Slepian'' basis functions are found for the
preimage which eventually leads to the calculation of an svd of the
restricted operator~$\cP\cT$. This also provides us with basis
functions which are orthogonal in the image spaces of $\cP\cT$ as well
as $\cT$. The singular values of $\cP\cT$ are linked to the maximized
norm quotient, and are diagnostic of the numerical stability and
ill-posedness of the inverse problem. We presented an algorithm for
determining the Slepian functions and the corresponding svd. We showed
how a wavelet multi-scale regularization can be constructed for a
variety of different filter functions.  Two numerical examples yielded
promising results. Our paper is an abstract generalization and an
illumination of the fundamental mathematical principles underlying the
method introduced in
\cite{PlattnerSimons2017}. In particular, we show how complicated problems with coupled sources can be integrated into our conceptual framework.

Practical examples where data are only regionally available or where
the analysis is only of interest in a particular subdomain are
abundant. In addition, we are often confronted with the situation that
the function of interest cannot be measured directly but is only
available via the solution of an ill-posed inverse problem. The
combination of both challenges (regional analysis and ill-posed
inverse problem) occurs rather often. We are now in the position to
further investigate the various possibilities that Slepian functions
provide for such inverse problems.

\section{Acknowledgments}

VM is grateful for the kind hospitality of Princeton University and,
in particular, of his host FJS during VM's sabbatical. The discussions
during our joint time at Princeton motivated this paper. This work was
partially supported by the Deutsche Forschungsgemeinschaft via grants
MI~655/7-2 and MI~655/10-1 to VM, and by the U.S. National Aeronautics
and Space Administration via grant NNX14AM29G to FJS and Alain
Plattner, and by the U.S. National Science Foundation via grant
EAR-1550389 to FJS and Alain Plattner.

\bibliographystyle{amsplain} 
\bibliography{Literatur}

\end{document}